\documentclass[ijoo,nonblindrev]{informs-ijoo}

\OneAndAHalfSpacedXI

\usepackage{nicefrac}
\usepackage{bm}
\usepackage{bbm}
\usepackage{enumitem}
\usepackage{comment}
\newcommand{\blank}{\kern.15em\cdot\kern.15em} 

\newcommand{\reals}{\mathbb{R}}

\newcommand{\nats}{\mathbb{N}}
\DeclareMathOperator{\bcfcns}{\mathcal{C}} 

\newcommand{\tsum}{\textstyle{\sum}}

\newcommand{\tto}{\rightrightarrows}

\DeclareMathAlphabet{\mathsl}{T1}{cmr}{m}{sl}
\newcommand{\drv}{\mathsl{d}}

\newcommand{\defeq}{\coloneq}

\DeclareMathOperator*{\minimize}{minimize}
\DeclareMathOperator*{\maximize}{maximize}
\let\argmin\relax
\DeclareMathOperator*{\argmin}{argmin}
\let\argmax\relax
\DeclareMathOperator*{\argmax}{argmax}

\newcommand{\subjectTo}{\text{subject to}}
\newcommand{\where}{\text{where}}

\newcommand{\rxi}{\bm{\xi}}
\newcommand{\romega}{\bm{\omega}}
\newcommand{\rp}{\bm{p}}
\newcommand{\rT}{\bm{T}}

\newcommand{\Expt}{\mathbb{E}} 
\newcommand{\Prob}{\mathbb{P}} 
\newcommand{\Qrob}{\mathbb{Q}} 
\newcommand{\pmProb}{\mathbbm{1}} 
\newcommand{\given}{\kern.1em\vert\kern.1em}
\newcommand{\givenn}{\kern-.15em\bigm\vert\kern-.15em}

\DeclareMathOperator{\Indicator}{\mathbbm{1}}

\newcommand{\stageC}{\varphi} 
\newcommand{\ctgF}{v} 
\newcommand{\ExptCtgF}{V} 
\newcommand{\oExptCtgF}{\widetilde{V}} 
\newcommand{\oCtgF}{\widetilde{v}} 
\newcommand{\SDP}{\mathrm{S}}
\newcommand{\MPC}{\mathrm{M}}
\newcommand{\DRO}{\mathrm{R}}
\newcommand{\DOO}{\mathrm{O}}

\newcommand{\storageCost}{C}
\newcommand{\storageCostDrv}{c}
\newcommand{\cPrice}{p^{\star}}
\newcommand{\cPriceSDP}{p_{\SDP}^{\star}}
\newcommand{\cPriceMPC}{p_{\MPC}^{\star}}

\usepackage{float} 

\usepackage{mathtools} 
\newcommand{\coloneq}{\coloneqq}
\let\subseteq\relax
\newcommand{\subseteq}{\subset}

\usepackage{hyperref}

\usepackage{microtype}

\usepackage{natbib}
 \bibpunct[, ]{(}{)}{,}{a}{}{,}%
 \def\bibfont{\small}%
\TheoremsNumberedThrough     
\ECRepeatTheorems

\EquationsNumberedThrough    

\MANUSCRIPTNO{IJO-2024-03-0040.R2}

\begin{document}


\RUNAUTHOR{Keehan, Philpott, and Anderson}

\RUNTITLE{Out-of-Sample Performance of SDP and MPC}

\TITLE{On the Out-of-Sample Performance of\\Stochastic Dynamic Programming\\and Model Predictive Control}

\ARTICLEAUTHORS{%
\AUTHOR{Dominic S. T. Keehan}
\AFF{University of Auckland, \EMAIL{dkee331@aucklanduni.ac.nz}} 
\AUTHOR{Andrew B. Philpott}
\AFF{University of Auckland, \EMAIL{a.philpott@auckland.ac.nz}}
\AUTHOR{Edward J. Anderson}
\AFF{Imperial College London, \EMAIL{e.anderson@imperial.ac.uk}}
} 

\ABSTRACT{
Sample average approximation--based stochastic dynamic programming (SDP) and model predictive control (MPC) are two different methods for approaching multistage stochastic optimization. In this paper we investigate the conditions under which SDP may be outperformed by MPC. We show that, depending on the presence of concavity or convexity, MPC can be interpreted as solving a mean-constrained distributionally ambiguous version of the problem that is solved by SDP. This furnishes performance guarantees when the true mean is known and provides intuition for why MPC performs better in some applications and worse in others. We then study a multistage stochastic optimization problem that is representative of the type for which MPC may be the better choice. We find that this can indeed be the case when the probability distribution of the underlying random variable is skewed or has enough weight in the right-hand tail.}

\KEYWORDS{Sample average approximation, stochastic dynamic programming, model predictive control, distributionally robust optimization.} 

\maketitle

%


\section{Introduction}
In practice multistage stochastic optimization problems often have to be solved without explicit knowledge of the probability distributions involved. Although one can create scenario-tree approximations of such problems based on samples of the random variables in each stage (termed \emph{sample average approximation} or SAA), the number of samples required to solve the true problem to a specified accuracy grows exponentially with the number of stages \citep{shapiro-nemirovski2005, shapiro2006} and solving the resulting problem is computationally expensive \citep{dyer-stougie2006}. It follows that for large problems SAA may only be practically applied when the number of samples in each stage is small. We are interested in the performance of different solution approaches in this small-sample regime.

Multistage stochastic optimization becomes easier when the random variables are stage-wise independent or follow a Markov process and the problem can be formulated as a stochastic optimal control problem. In principle, such problems are amenable to solution via stochastic dynamic programming (SDP), as long as the state dimension is not too high \citep{dynamic-programming}. However, there are alternative solution approaches. Often, model predictive control (MPC) is used. MPC fixes the random variables involved in the  optimization problem to deterministic forecasts, giving a cost-to-go function that is computationally inexpensive to evaluate \citep{bertsekas2005}. In this way MPC can handle high state dimensions, nonlinear constraints, and a large number of stages. Although the optimal solutions obtained from SDP and MPC coincide for certain quadratic problems \citep{theil1957, ziemba1971}, this is the exception rather than the rule.

In our paper both SDP and MPC are sample based. We use SAA to construct approximations of the true probability distributions, assigning {equal probability to each of} $N$ random samples. Also, we use a simple form of MPC that forecasts the random variables as the average of these $N$ samples. Compared to SDP, in certain applications MPC performs poorly out-of-sample \citep{pacaud-et-al-2024-mpc-sddp}, and in other applications MPC performs well out-of-sample \citep{martin-PhD-2021}. This is surprising as MPC is tantamount to ignoring randomness.

The intention of this paper is to advance our understanding of SDP and MPC applied to stochastic optimal control problems, by comparing the out-of-sample performance of each method when $N$ is small. Our study is motivated by the question:
\medskip
\begin{center}
    \emph{Under what conditions is SDP outperformed by MPC out of sample?}
\end{center}
\medskip
To answer this question we first look through the lens of distributional ambiguity. In data-driven stochastic optimization, when $N$ is small, considering an ambiguity set of probability distributions informed by the $N$ samples can improve out-of-sample performance \citep{anderson-philpott-2022, gotoh-etal-2023-dro-doo}. Moment constraints are a classical choice for constructing these ambiguity sets \citep{scarf1958, dupacova-1966-minimax, delage-ye-2010-moment-dro}. Owing to Jensen's inequality, depending on the convexity or concavity of the cost-to-go function in the random variables, we show that sample-based MPC can be interpreted as solving a mean-constrained distributionally ambiguous version of the problem that is solved by sample-based SDP.

The distributionally ambiguous interpretations we provide are {not} guarantees that one method will outperform the other, but they provide an intuition for the types of problems for which one method may be a better choice than the other. To gain a deeper understanding of the conditions under which SDP can be outperformed by MPC out of sample, we study a specific revenue optimization problem that determines how much product from an existing inventory to sell at a given price offered by the market. Unsold product may be stored as inventory for later sales at a different realisation of the random price, subject to deterministic holding costs. This problem is simple enough to admit the derivation of a closed-form optimal policy, but complex enough to capture critical aspects of the differences in out-of-sample performance between the methods.

Given the revenue optimization problem and some ground-truth probability distribution for the random price, we derive an optimal SDP policy for any $N$ fixed samples and hence can evaluate its out-of-sample performance under the true distribution. Similarly, we derive an optimal MPC policy and can evaluate its performance under the true distribution. These two values enable us to understand the sensitivities of each method to the specific values of the $N$ samples used in the distributional approximations. Furthermore, the expectation of these two values over the sampling distribution of the $N$ samples quantifies the average performance of each method.

In Section~\ref{section:stochastic-optimal-control} we formulate a general stochastic optimal control problem as well as its data-driven approximations via sample-based SDP and MPC. In Section~\ref{section:MPC-is-DRO} we classify the Bellman operator associated with the MPC problem as concavity or convexity preserving, and show that in the former case MPC can be interpreted as solving a mean-constrained distributionally robust version of the problem that is solved by SDP. In the latter case this becomes a distributionally optimistic version. We then consider a number of specific examples where the MPC Bellman operator is concavity or convexity preserving. In Section~\ref{section:MPC-guarantee-and-performance} we show that the distributionally ambiguous interpretations developed in the previous section can be used to derive performance guarantees for MPC when the true mean is known. We also derive a result for comparing the out-of-sample performance of different policies for stochastic optimal control problems. In Section~\ref{section:multistage-revenue-optimisation} we introduce our multistage revenue optimization problem and establish a closed-form expression for its optimal solution. We then compare the out-of-sample performance of SDP and MPC and provide a condition on the samples which ensures that MPC performs at least as well as SDP. In the remainder of Section~\ref{section:multistage-revenue-optimisation} we report on some examples. Firstly, we show that when the underlying distribution of the random price is exponential, the expected out-of-sample performance improvement from using MPC instead of SDP becomes arbitrarily large as the discount factor approaches $1$. Lastly, we present a range of numerical experiments that support the observations of the previous sections. The paper concludes with a discussion in Section~\ref{section:SDP-vs-MPC-discussion}. 

All of the proofs of the results in this
paper are deferred to the appendices.

\section{Stochastic Optimal Control}
\label{section:stochastic-optimal-control}
Before introducing the control problem, we give some notation. An element $x\in\mathcal{X} \subset \reals^n$ denotes the state vector, and an element $\xi\in\Xi \subset \reals^m$ denotes the outcome of a random noise vector $\rxi$ distributed according to $\Prob$. (Boldface distinguishes random vectors from their outcomes.) Superscripts denote indexing by stage, but also, where clear from context, raising to a power.

For initial state $x^{1}\in\mathcal{X}$ we study the infinite-horizon stochastic optimal control problem
\begin{alignat}{2}\label{problem:soc-sum}
& \kern 0.04em \minimize_{y^1, y^2, \ldots} && \quad \Expt_{\Prob^{\infty}}\Biggl[\,\sum_{t=1}^{\infty}\beta^{t-1}\stageC(x^t,x^{t+1},\rxi^t) \Biggr]
\tag{SOC}\\
& \subjectTo &&
\quad x^{t+1} = y^{t}(x^t,\rxi^t) \in \mathcal{Y}(x^t,\rxi^t) \quad  \forall t\in \nats, \notag
\end{alignat}
where $\{y^t:\mathcal{X}\times\Xi \to\mathcal{X},\,t\in \nats\}$ is a \emph{policy} that at each stage $t$ provides a \emph{decision rule} $y^t$ for the next state $x^{t+1}$ given the current state $x^t$ and the realisation of the random vector $\rxi^t$. At each stage the mapping $\xi \mapsto y^t(x,\xi)$ is required to be measurable for every $x\in\mathcal{X}$. Here $\beta \in (0,1)$ is a discount factor (with $\beta^t$ its $t$\textsuperscript{th} power), $\stageC:\mathcal{X}\times\mathcal{X}\times\Xi\to\reals$ is a cost function, and $\mathcal{Y}:\mathcal{X}\times\Xi\tto \mathcal{X}$ is a set-valued mapping enforcing state-transition constraints. The infinite-product distribution $\Prob^{\infty}\defeq \Prob \times \Prob\times\cdots$ defines a joint distribution for the random vectors $\rxi^1,\rxi^2,\ldots\,$.

In this formulation we have assumed that the dynamics of the state evolution is such that the control is exercised after the noise is known. For that reason we can take the control as simply the choice of state at the next period. The noise has an effect on the possible choice of next state, as well as on the cost function.

In general, (\ref{problem:soc-sum}) may not be well defined. To keep our analysis simple, we make the following assumption.
\begin{assumption}\label{assumption:bounded-continuous}
\hfill \begin{enumerate}[label={\normalfont(\roman*)}, ref={\roman*}]
\item\label{assumption:bounded-continuous-i} The set $\mathcal{X}$ is closed and convex, and the set-valued mapping $\mathcal{Y}:\mathcal{X}\times\Xi\tto\mathcal{X}$ is nonempty- and compact-valued, and continuous.
\item\label{assumption:bounded-continuous-ii} The set $\Xi$ is closed and convex, and the random vectors $\rxi^1,\rxi^2,\ldots$ are independent and identically distributed according to $\Prob$.
\item\label{assumption:bounded-continuous-iii} The function $\stageC:\mathcal{X}\times\mathcal{X}\times\Xi\to\reals$ is bounded and continuous.
\end{enumerate}
\end{assumption}
Note that Assumption~\ref{assumption:bounded-continuous}(\ref{assumption:bounded-continuous-ii}) does not preclude modelling stagewise-dependent randomness, since auxiliary states that track the realisations of random vectors in previous stages can be created, as we do in Example~\ref{example:AR1-objective-uncertainty} below. Moreover, Assumption~\ref{assumption:bounded-continuous}(\ref{assumption:bounded-continuous-iii}) can often be relaxed, although the mathematical tools required become more complex.

For simplicity we do not vary $\Prob$, $\stageC$, or $\mathcal{Y}$ between stages. However, our results can be generalised to such a setting using a periodic formulation of (\ref{problem:soc-sum}), as in \citep{shapiro-periodic-consistency}. Also, in this paper we consider only infinite-horizon stochastic optimal control problems, although all of our results have finite-horizon counterparts. For a comparison of the finite- and infinite-horizon cases, see \citep[Chapters~3 and 4]{bertsekas1996stochastic}.

The problem (\ref{problem:soc-sum}) is an infinite-dimensional optimization problem. A closely related problem is that of finding a function $\ctgF:\mathcal{X}\times\Xi\to\reals$ which satisfies the functional equation
\begin{equation}
     \ctgF(x,\xi)= \inf_{y\in\mathcal{Y}(x,\xi)}\Bigl\{\stageC(x,y,\xi)+\beta \Expt_{\Prob} \bigl[\ctgF(y,\rxi)\bigr]\Bigr\} \quad \forall (x,\xi)\in\mathcal{X}\times\Xi.
\label{equation:value-function}
\end{equation}
Under Assumption~\ref{assumption:bounded-continuous} the equation (\ref{equation:value-function}) has a unique bounded and continuous solution $\ctgF:\mathcal{X}\times\Xi\to\reals$. Moreover, the value $\Expt_{\Prob}\bigl[\ctgF(x^1,\rxi)\bigr]$ is equal to the optimal value of (\ref{problem:soc-sum}), and there is a function  $y:\mathcal{X}\times\Xi\to\mathcal{X}$ with $\xi\mapsto y(x,\xi)$ measurable for every $x\in\mathcal{X}$ satisfying 
\begin{equation*}
y(x,\xi)\in\argmin_{y\in\mathcal{Y}(x,\xi)}\Bigl\{\stageC(x,y,\xi)+\beta\Expt_{\Prob}\bigl[\ctgF(y,\rxi)\bigr]\Bigr\} \quad \forall (x,\xi)\in\mathcal{X}\times\Xi.
\end{equation*}
This defines a decision rule for each stage and the resulting policy solves (\ref{problem:soc-sum}). For a discussion of all of these well-known results, see, e.g., \citep[Chapter~9]{stokey-lucas:recursive-methods} or \citep[Chapter~4]{bertsekas1996stochastic}. We refer to any function $y:\mathcal{X}\times\Xi\to\mathcal{X}$ with $\xi\mapsto y(x,\xi)$ measurable and $y(x,\xi)\in\mathcal{Y}(x,\xi)$ for every $(x,\xi)\in\mathcal{X}\times\Xi$ as an \emph{admissible policy} for (\ref{problem:soc-sum}).

\subsection{Data-Driven Stochastic Optimal Control}
To approximate (\ref{problem:soc-sum}) in a data-driven manner, given $N$ samples $\xi_1,\ldots,\xi_N\in\Xi$ we form the empirical distribution $\Prob_N$. The infinite-product distribution $\Prob_N^{\infty} = \Prob_N\times \Prob_N\times\cdots$ then defines an approximation of the joint distribution for the random vectors $\rxi^1,\rxi^2,\ldots\,$, which yields the problem
\begin{alignat}{2}\label{problem:saa-sum}
& \kern 0.04em \minimize_{y^1, y^2, \ldots} && \quad  \Expt_{\Prob_N^{\infty}}\Biggl[\,\sum_{t=1}^{\infty} \beta^{t-1} \stageC(x^t,x^{t+1},\rxi^t) \Biggr]\\
& \subjectTo &&
\quad x^{t+1} = y^{t}(x^t,\rxi^t) \in \mathcal{Y}(x^t,\rxi^t) \quad  \forall t\in \nats.\notag
\end{alignat}
As for (\ref{equation:value-function}), this is related to the problem of finding a function ${\ctgF_{\SDP}:\mathcal{X}\times\Xi\to\reals}$ which satisfies the sample-based SDP functional equation
\begin{equation}
     \ctgF_{\SDP}(x,\xi)= \inf_{y\in\mathcal{Y}(x,\xi)}\Bigl\{\stageC(x,y,\xi)+\beta \Expt_{\Prob_N} \bigl[\ctgF_{\SDP}(y,\rxi)\bigr]\Bigr\} \quad \forall (x,\xi)\in\mathcal{X}\times\Xi.
\label{equation:sdp-value-function}
\end{equation}
Under Assumption~\ref{assumption:bounded-continuous} the same remarks regarding uniqueness, boundedness, continuity, and measurability made with reference to (\ref{equation:value-function}) apply to (\ref{equation:sdp-value-function}), and by solving (\ref{equation:sdp-value-function}) a sample-based SDP policy $y_{\SDP}:\mathcal{X}\times\Xi\to\mathcal{X}$ satisfying 
\begin{equation*}
y_{\SDP}(x,\xi)\in\argmin_{y\in\mathcal{Y}(x,\xi)}\Bigl\{\stageC(x,y,\xi)+\beta\Expt_{\Prob_N}\bigl[\ctgF_{\SDP}(y,\rxi)\bigr]\Bigr\} \quad \forall (x,\xi)\in\mathcal{X}\times\Xi
\end{equation*}
can be obtained.

MPC proceeds differently and replaces the random vectors in (\ref{problem:soc-sum}) with a deterministic forecast. Given $N$ samples $\xi_1,\ldots,\xi_N\in\Xi$, we use the sample average $\mu_N \defeq \frac{1}{N}\sum_{i=1}^N \xi_i \in \Xi$ as the forecast. This leads to the problem of finding a function $\ctgF_{\MPC}:\mathcal{X}\times\Xi\to\reals$ which satisfies the sample-based MPC functional equation
\begin{equation}
     \ctgF_{\MPC}(x,\xi)= \inf_{y\in\mathcal{Y}(x,\xi)}\Bigl\{\stageC(x,y,\xi)+\beta \ctgF_{\MPC}(y,\mu_N )\Bigr\} \quad \forall (x,\xi)\in\mathcal{X}\times\Xi.
\label{equation:mpc-value-function}
\end{equation}
Observe that (\ref{equation:mpc-value-function}) is equivalent to (\ref{equation:value-function}) with the distribution $\pmProb_{\mu_N}$ used instead of $\Prob$. (Here $\pmProb_{\xi}$  denotes the point-mass distribution that assigns probability $1$ to $\xi\in\Xi$.) Under Assumption~\ref{assumption:bounded-continuous}, the same remarks regarding uniqueness, boundedness, continuity, and measurability made with reference to (\ref{equation:value-function}) apply to (\ref{equation:mpc-value-function}), and by solving (\ref{equation:mpc-value-function}) a sample-based MPC policy ${y_{\MPC}:\mathcal{X}\times\Xi\to\mathcal{X}}$ satisfying 
\begin{equation*}
y_{\MPC}(x,\xi)\in\argmin_{y\in\mathcal{Y}(x,\xi)}\Bigl\{\stageC(x,y,\xi)+\beta\ctgF_{\MPC}(y,\mu_N )\Bigr\} \quad \forall (x,\xi)\in\mathcal{X}\times\Xi
\end{equation*}
can be obtained. We are interested in comparing the performance of the sample-based SDP and MPC policies when applied out of sample to the true problem (\ref{problem:soc-sum}).

\section{Connections With Distributionally Ambiguous Optimization}\label{section:MPC-is-DRO}
The SDP policy obtained by solving (\ref{equation:sdp-value-function}) depends critically on the values of the $N$ samples used. Instead, by hedging against the worst-case distributions in an \emph{ambiguity set} $\mathcal{P}$ informed by the $N$ samples, \emph{distributionally robust optimization} (DRO) provides a means to limit out-of-sample disappointment. The infinite product of these ambiguity sets $\mathcal{P}^\infty \defeq \mathcal{P} \times \mathcal{P} \times \cdots$ defines a set of joint distributions for the random vectors $\rxi^1,\rxi^2,\ldots\,$, which yields the distributionally robust problem
\begin{alignat}{2}\label{problem:dro-sum}
& \kern 0.04em \minimize_{y^1, y^2, \ldots} && \quad \sup_{\Qrob\in\mathcal{P}^{\infty}} \Expt_{\Qrob} \Biggl[\,\sum_{t=1}^{\infty} \beta^{t-1} \stageC(x^t,x^{t+1},\rxi^t) \Biggr]\\
& \subjectTo &&
\quad x^{t+1} = y^{t}(x^t,\rxi^t) \in \mathcal{Y}(x^t,\rxi^t) \quad  \forall t\in \nats.\notag
\end{alignat}

The choice of $\mathcal{P}^\infty$ in (\ref{problem:dro-sum}) requires care. By using a product set, under the assumption of stage-wise independence, $\mathcal{P}^\infty$ satisfies the ``rectangularity'' property  (see \citep[Assumption~2.1]{iyengar-2005-robust-dp} and \citep{shapiro-2016-rectangularity} for details), whereby (\ref{problem:dro-sum}) can be solved by finding a function $\ctgF_{\DRO}:\mathcal{X}\times\Xi\to\reals$ which satisfies the DRO functional equation
\begin{equation}
     \ctgF_{\DRO}(x,\xi)= \inf_{y\in\mathcal{Y}(x,\xi)}\Bigl\{\stageC(x,y,\xi)+\beta 
     \sup_{\Qrob\in\mathcal{P}} \Expt_{\Qrob} \bigl[\ctgF_{\DRO}(y,\rxi)\bigr]\Bigr\} \quad \forall (x,\xi)\in\mathcal{X}\times\Xi.
\label{equation:dro-value-function}
\end{equation}
Note that the ambiguity set in (\ref{equation:dro-value-function}) is $\mathcal{P}$ so the supremum is over a set of distributions for $\rxi$, rather than a set of joint distributions for $\rxi^1,\rxi^2,\ldots$ as in (\ref{problem:dro-sum}).

To demonstrate a connection between DRO and MPC, we now analyse the functional operator associated with (\ref{equation:mpc-value-function}). Let $\bcfcns(\mathcal{X}\times\Xi)$ denote the set of bounded and continuous real-valued functions on $\mathcal{X}\times\Xi$. This is a Banach space under the sup norm $\lVert\blank\rVert_{\infty}$. Under Assumption~\ref{assumption:bounded-continuous}, the equation (\ref{equation:mpc-value-function}) features an MPC Bellman operator $B_{\MPC}:\bcfcns(\mathcal{X}\times\Xi)\to\bcfcns(\mathcal{X}\times\Xi)$ which for $f \in \bcfcns(\mathcal{X}\times\Xi)$ at $(x,\xi)\in\mathcal{X}\times\Xi$ has value
\begin{equation*}\label{equation:mpc-bellman-operator}
    B_{\MPC}(f)(x,\xi) \defeq \inf_{y\in\mathcal{Y}(x,{\xi})}\Bigl\{\stageC(x,y,{\xi})+\beta f(y,\mu_N )\Bigr\}.
\end{equation*} 
The operator $B_{\MPC}$ can be defined for any forecast in place of the sample average $\mu_N$ which appears here. For a proof of the fact that $B_{\MPC}$ maps $\bcfcns(\mathcal{X}\times\Xi)$ to itself, see, e.g., \citep[Theorem~9.6]{stokey-lucas:recursive-methods}.

Recall that $B_{\MPC}$ is a contraction mapping under the sup norm, since
\begin{equation*}
    \lVert B_{\MPC}(f) - B_{\MPC}(g) \rVert_{\infty} \leq \beta \lVert f - g \rVert_{\infty} \quad \forall f,g\in\bcfcns(\mathcal{X}\times\Xi)
\end{equation*}
with $\beta\in(0,1)$. Banach's fixed-point theorem then shows that the equation $f = B_{\MPC}(f)$ has a unique solution in $(\bcfcns(\mathcal{X}\times\Xi),\lVert\blank\rVert_{\infty})$, and this solution solves the MPC functional equation (\ref{equation:mpc-value-function}). Further, if there exists a closed subset $\mathcal{F}\subseteq\bcfcns(\mathcal{X}\times\Xi)$ which maps to itself under $B_{\MPC}$ and ${f\in\bcfcns(\mathcal{X}\times\Xi)}$ solves $f = B_{\MPC}(f)$, then it holds that $f\in\mathcal{F}$ \citep[Theorem~3.2, Corollary~1]{stokey-lucas:recursive-methods}. This motivates the following definition.
\begin{definition}
The MPC Bellman operator $B_{\MPC}$ is \emph{concavity preserving} (respectively \emph{convexity preserving}) if there exists a closed subset $\mathcal{F}\subseteq\bcfcns(\mathcal{X}\times\Xi)$  where the mapping $\xi\mapsto f(x,\xi)$ is concave (convex) for every $x\in\mathcal{X}$ and $f\in \mathcal{F}$, and which maps to itself under $B_{\MPC}$.
\end{definition}

As we will show, problems with concavity-preserving MPC Bellman operators give rise to sample-based MPC policies that satisfy a DRO interpretation. A sufficient condition for this concavity preservation is provided by the following result.
\begin{proposition}
\label{proposition:concave-costs}
Let Assumption~\ref{assumption:bounded-continuous} hold. Suppose that the set-valued mapping $\mathcal{Y}$ does not depend on $\xi$ and that the mapping $\xi\mapsto\stageC(x,y,\xi)$ is concave for 
every $x, y \in \mathcal{X}$. Then the MPC Bellman operator $B_{\MPC}$ is concavity preserving.
\end{proposition}
\noindent In some circumstances we can establish this key property of concavity preservation without requiring all the conditions of Proposition~\ref{proposition:concave-costs}. Some special cases will be discussed in examples below.

When the MPC Bellman operator is concavity preserving, we have the following DRO interpretation.
\begin{theorem}
Let Assumption~\ref{assumption:bounded-continuous} hold. For an ambiguity set $\mathcal{P}$, suppose that $\Expt_{\Qrob}[\rxi] = \mu_N$ for every $\Qrob\in\mathcal{P}$ and that $\pmProb_{\mu_N} \in \mathcal{P}$. If the MPC Bellman operator $B_{\MPC}$ is concavity preserving, then the DRO functional equation (\ref{equation:dro-value-function}) has a solution which is the same as the solution to the MPC functional equation (\ref{equation:mpc-value-function}).
\label{theorem:mpc-is-dro}
\end{theorem}

In contrast to DRO, \emph{distributionally optimistic optimization} (DOO) considers only the best-case distributions in an ambiguity set $\mathcal{P}$. Using the rectangular ambiguity set $\mathcal{P}^{\infty}$, a distributionally optimistic version of (\ref{problem:dro-sum}) can be formed. This leads to the problem of finding a function ${\ctgF_{\DOO}:\mathcal{X}\times\Xi\to\reals}$ which satisfies the DOO functional equation
\begin{equation}
     \ctgF_{\DOO}(x,\xi)= \inf_{y\in\mathcal{Y}(x,\xi)}\Bigl\{\stageC(x,y,\xi)+\beta 
     \inf_{\Qrob\in\mathcal{P}} \Expt_{\Qrob} \bigl[\ctgF_{\DOO}(y,\rxi)\bigr]\Bigr\} \quad \forall (x,\xi)\in\mathcal{X}\times\Xi.
\label{equation:doo-value-function}
\end{equation}

The essential characteristic of the connection between DRO and MPC is the concavity-preserving nature of the operator $B_{\MPC}$. If $B_{\MPC}$ is instead convexity preserving, the following converse result to Theorem~\ref{theorem:mpc-is-dro} holds. 
\begin{corollary}
Suppose that the conditions in the statement of Theorem~\ref{theorem:mpc-is-dro} hold, but that the MPC Bellman operator $B_{\MPC}$ is convexity preserving rather than concavity preserving. Then the DOO functional equation (\ref{equation:doo-value-function}) has a solution which is the same as the solution to the MPC functional equation (\ref{equation:mpc-value-function}).
\label{corollary:mpc-is-doo}
\end{corollary}

We now give some examples of problems with concavity- and convexity-preserving MPC Bellman operators that differ from those characterised by Proposition~\ref{proposition:concave-costs}.

\begin{example}\label{example:AR1-objective-uncertainty}
\citet{downward_et_al_2020_objective_uncertainty} study a class of problems where stagewise-dependent uncertainty appears in the objective function \citep[see also][Chapter~4]{dowson-PhD}. These include resource allocation problems under price uncertainty, such as hydroelectric bidding and agricultural production planning. We present an example from this class that minimizes the production cost of a commodity sold on a spot market. Suppose at time $t$ the spot price $p^t$ is an autoregressive process defined by $\rp^{t}=\alpha p^{t-1} +\rxi$. Further, suppose the commodity has per-period production $w$ and marginal per-period storage cost $c$. Here the state is determined by the inventory $x\in\mathcal{X}$, and now also by the previous spot price $q\in\mathcal{Q}$. As before the decision variables are the next state: an inventory-price pair $(y,p)$, but in this case $p$ is the current price, already determined by $q$ and $\xi$. The functional equation is
\begin{equation*}
\ctgF(x,q,\xi)=\inf_{(y,p) \in \mathcal{Y}(x,q,\xi)}\Bigl\{- p(x+w-y) + cy +\beta \Expt\bigl[\ctgF(y,p,\rxi)\bigr]\Bigr\} \quad \forall (x,q,\xi)\in\mathcal{X}\times\mathcal{Q}\times\Xi
\end{equation*}
where
\begin{equation*}
\mathcal{Y}(x,q,\xi) = \bigl\{(y,p)\in\mathcal{X}\times\mathcal{Q} : y \leq x+w,\;  p=\alpha q+\xi\bigr\}.
\end{equation*}
To see that the MPC Bellman operator is concavity preserving, consider the set
\begin{equation*}
    \mathcal{F} = \bigl\{f\in\bcfcns(\mathcal{X}\times\mathcal{Q}\times\Xi) : (q,\xi)\mapsto f(x,q,\xi)~\text{is concave for every}~ x\in\mathcal{X}\bigr\}.
\end{equation*}
For every $f \in \mathcal{F}$ the mapping $p\mapsto -p(x+w-y)+cy+\beta f(y,p,\Expt[\rxi])$ is concave. Due to the equality constraint $p=\alpha q + \xi$, it can be shown that the infimum operation preserves concavity in $(q,\xi)$ through a similar argument to that in the proof of Proposition~\ref{proposition:concave-costs}. Thus, sample-based MPC satisfies the DRO interpretation of Theorem~\ref{theorem:mpc-is-dro}. \Halmos
\end{example}




\begin{example}\label{example:convexity-preserving}
The \emph{hydrothermal scheduling problem} \citep{pereira1991multi} minimizes the thermal energy cost of meeting energy demand in a system with hydroelectric plants supplied by reservoir storage with random inflows. Consider a thermal plant with generation $t\in\mathcal{T}$, constant per-period demand $d$, and a hydro plant that generates $h\in\mathcal{H}$ from a reservoir with storage $x\in\mathcal{X}$ and random inflow $\rxi$. All quantities are measured in terms of energy. The functional equation is
\begin{equation*}
\ctgF(x,\xi)=\inf_{(y,h,t) \in \mathcal{Y}(x,\xi)}
\Bigl\{t+\beta \Expt\bigl[\ctgF(y,\rxi)\bigr]\Bigr\} \quad \forall (x,\xi)\in\mathcal{X}\times\Xi
\end{equation*}
where
\begin{equation*}
\mathcal{Y}(x,\xi) = \bigl\{(y,h,t)\in\mathcal{X}\times\mathcal{H}\times\mathcal{T}:h+t \geq d,\;h\leq x,\; y = x - h + \xi\bigr\}.
\end{equation*}
Here $h$ and $t$ are viewed as dummy state variables in our framework. To see that the MPC Bellman operator is convexity preserving, consider the set
\begin{equation*}
    \mathcal{F} = \bigl\{f\in\bcfcns(\mathcal{X}\times\Xi):(x,\xi)\mapsto f(x,\xi)~\text{is convex}\bigr\}.
\end{equation*}
For every $f \in \mathcal{F}$ the mapping $(y,h,t)\mapsto t+\beta f(y,\Expt[\rxi])$ is convex. Due to the linearity of the constraints, it can be shown that the infimum operation preserves convexity in $(x,\xi)$ through a similar argument to that in the proof of Proposition~\ref{proposition:concave-costs} (although with the inequalities reversed). Thus, sample-based MPC satisfies the DOO interpretation of Corollary~\ref{corollary:mpc-is-doo}. \Halmos
\end{example}





If the MPC Bellman operator is concavity preserving, Theorem~\ref{theorem:mpc-is-dro} shows that sample-based MPC can be interpreted as a distributionally robust version of sample-based SDP. Distributional robustness can sometimes improve out-of-sample performance due to the effects of shrinkage \citep{anderson-philpott-2022}. Conversely, if the MPC Bellman operator is convexity preserving, Corollary~\ref{corollary:mpc-is-doo} shows that sample-based MPC can be interpreted as a distributionally optimistic version of sample-based SDP. Distributional optimism can sometimes worsen out-of-sample performance due to increased solution variance \citep{gotoh-etal-2023-dro-doo}. This suggests that for small sample sizes, when the MPC Bellman operator is concavity preserving as in Example~\ref{example:AR1-objective-uncertainty}, MPC may be a better choice than SDP. On the other hand, when the MPC Bellman operator is convexity preserving as in Example~\ref{example:convexity-preserving}, SDP may be a better choice than MPC. In the following section we provide performance guarantees which support these observations.

\section{Out-of-Sample Performance}
\label{section:MPC-guarantee-and-performance}
Suppose that the distribution $\Prob$ has a finite mean $\Expt_{\Prob}[\rxi]$. For $N$ independent random samples $\rxi_1,\ldots,\rxi_N\sim\Prob$, the sample average $\frac{1}{N}\sum_{i=1}^N\rxi_i$ is usually a good estimate of the true mean $\Expt_{\Prob}[\rxi]$, as long as $N$ is sufficiently large. Indeed, the strong law of large numbers shows ${\frac{1}{N}\sum_{i=1}^N\rxi_i \to \Expt_{\Prob}[\rxi]}$ almost surely as $N\to\infty$. MPC can take advantage of this situation since the computational complexity of solving (\ref{equation:mpc-value-function}) does not depend on the sample size, while SDP cannot.

We now show that, when the true mean is known and used as the MPC forecast, the solution to the MPC functional equation~(\ref{equation:mpc-value-function}) provides a performance guarantee on the cost incurred when applying the resulting MPC policy to (\ref{problem:soc-sum}); that is, we will do better than predicted.
\begin{theorem}\label{theorem:mpc-dro-guarantee}
Let Assumption~\ref{assumption:bounded-continuous} hold and set $\mu_N=\mu \defeq \Expt_{\Prob}[\rxi]$ in the MPC functional equation~(\ref{equation:mpc-value-function}). If the resulting MPC Bellman operator $B_\MPC$ is concavity preserving, then
\begin{equation*}
    \Expt_{\Prob^\infty}\Biggl[\,\sum_{t=1}^{\infty}\beta^{t-1} \stageC(x^t,x^{t+1},\rxi^t)\Biggr] \leq \ctgF_{\MPC}(x^1,\mu) \quad \where~x^{t+1} = y_{\MPC}(x^t,\rxi^t) \; \forall t\in \nats.
\end{equation*}
That is, the value $\ctgF_{\MPC}(x^1,\mu)$ obtained by solving the MPC functional equation (\ref{equation:mpc-value-function}) is an upper bound on the cost incurred when applying the resulting policy to (\ref{problem:soc-sum}).
\end{theorem}

If the sample size is large then using MPC with the sample mean as the forecast is likely to give a result close to that with the true mean. This suggests that applying the sample-based MPC policy to (\ref{problem:soc-sum}) is unlikely to cause out-of-sample disappointment when the MPC Bellman operator is concavity preserving. 

We also have a lower bound on costs in the case where the MPC Bellman operator is convexity preserving, so applying MPC with the true mean used as the forecast will do worse than predicted.
\begin{corollary}
Suppose that the conditions in the statement of Theorem~\ref{theorem:mpc-dro-guarantee} hold, but that the MPC Bellman operator is convexity preserving rather than concavity preserving. Then
\begin{equation*}
    \Expt_{\Prob^\infty}\Biggl[\,\sum_{t=1}^{\infty}\beta^{t-1} \stageC(x^t,x^{t+1},\rxi^t)\Biggr] \geq \ctgF_{\MPC}(x^1,\mu) \quad \text{where}~x^{t+1} = y_{\MPC}(x^t,\rxi^t) \;\forall t\in \nats.
\end{equation*}
That is, the value $\ctgF_{\MPC}(x^1,\mu)$ obtained by solving the MPC functional equation (\ref{equation:mpc-value-function}) is a lower bound on the cost incurred when applying the resulting policy to (\ref{problem:soc-sum}).
\label{corollary:mpc-doo-dissapointment}
\end{corollary}
\noindent Corollary~\ref{corollary:mpc-doo-dissapointment} suggests that sample-based MPC is likely to cause out-of-sample disappointment when the MPC Bellman operator is convexity preserving.

In the rest of this section we derive a result for comparing the cost of applying different policies to (\ref{problem:soc-sum}). Given an admissible policy $y:\mathcal{X}\times\Xi\to\mathcal{X}$, with $\stageC$ continuous, the mapping ${\xi \mapsto \stageC(x,y(x,\xi),\xi)}$ is measurable for every $x\in\mathcal{X}$. Hence, the objective value of (\ref{problem:soc-sum}) under $y$ is well defined. Let $\oExptCtgF_{y}:\mathcal{X}\to\reals$ denote this value as a function of the initial state, so
\begin{equation}\label{equation:out-of-sample-sum}
    \oExptCtgF_{y}(x^1) \defeq \Expt_{\Prob^\infty}\Biggl[\,\sum_{t=1}^{\infty}\beta^{t-1} \stageC(x^t,x^{t+1},\rxi^t)\Biggr] \quad \text{where}~x^{t+1} = y(x^t,\rxi^t) \;  \forall t\in \nats.
\end{equation}
If $y$ is constructed from samples, then $\oExptCtgF_y(x^1)$ is the \emph{out-of-sample cost} of $y$.

To study out-of-sample performance, we do not require the boundedness condition of Assumption~\ref{assumption:bounded-continuous}(\ref{assumption:bounded-continuous-iii}). Instead, we make the following integrability assumption.
\begin{assumption}\label{assumption:unbounded-noise}
Assumption~\ref{assumption:bounded-continuous}(\ref{assumption:bounded-continuous-i}) and (\ref{assumption:bounded-continuous-ii}) hold, and
\begin{enumerate}[label={\normalfont(\roman*)}, ref={\roman*}]
\setcounter{enumi}{2}
\item\label{assumption:unbounded-noise-iii} The function $\stageC:\mathcal{X}\times\mathcal{X}\times\Xi\to\reals$ is continuous, and there exists a positive-valued random variable $l(\rxi)$ with $\Expt_{\Prob}\bigl[l(\rxi)\bigr]< \infty$, such that for every $x,y\in\mathcal{X}$ and $\Prob$-almost every $\xi\in\Xi$ it holds that $l(\xi) \geq \bigl\lvert\stageC(x,y,\xi)\bigr\rvert$.
\end{enumerate}
\end{assumption}

Under Assumption~\ref{assumption:unbounded-noise}(\ref{assumption:unbounded-noise-iii}), the $t$\textsuperscript{th} term in the sum in (\ref{equation:out-of-sample-sum}) is bounded by the random variable $l(\rxi^t)$ which has a finite expectation. With $\beta \in (0,1)$, it follows that $\oExptCtgF_{y}$ is bounded. Having defined $\oExptCtgF_{y}$ as a function of the initial state and shown that it is bounded, it satisfies the functional equation
\begin{equation}
\label{equation:out-of-sample-functional-equation}
\oExptCtgF_{y}(x)=\Expt_{\Prob}\Bigl[\stageC\bigl(x,y(x,\rxi),\rxi\bigr)
+\beta \oExptCtgF_{y}\bigl(y(x,\rxi)\bigr)\Bigr] \quad \forall x\in\mathcal{X}.
\end{equation}

Our approach to compare two different policies is to consider starting with one policy and then switching to the other policy after a certain number of stages. To this end, we make the following definition.
\begin{definition}\label{definition:switch-policies-value-function}
For two admissible policies $y,y^{\prime}:\mathcal{X}\times\Xi\to\mathcal{X}$, define $\oExptCtgF_{y,y^{\prime}}:\mathcal{X}\to\reals$ by
\begin{equation*}
\oExptCtgF_{y,y^{\prime}}(x) \defeq \Expt_{\Prob}\Bigl[ \stageC\bigl(x,y(x,\rxi),\rxi\,\bigr)
+\beta\oExptCtgF_{y^{\prime}}\bigl(y(x,\rxi)\bigr)\Bigr].
\end{equation*}
\end{definition}
The value $\oExptCtgF_{y,y^{\prime}}(x)$ is the objective value of (\ref{problem:soc-sum}) when starting from initial state $x$ if policy $y$ is followed in the first stage and policy $y^{\prime}$ is followed thereafter. It is clear that $\oExptCtgF_{y,y^{\prime}}$ is well defined and bounded in the same way that $\oExptCtgF_{y}$ and $\oExptCtgF_{y^{\prime}}$ are.

\begin{theorem}\label{theorem:policy-switch-improvement}
Let Assumption~\ref{assumption:unbounded-noise} hold. For two admissible policies $y,y^{\prime}:\mathcal{X}\times\Xi\to\mathcal{X}$, if
\begin{equation*}
     \oExptCtgF_{y^{\prime}}(x)\leq\oExptCtgF_{y,y^{\prime}}(x) \quad \forall x \in \mathcal{X},
\end{equation*}
then $\oExptCtgF_{y^{\prime}}(x)\leq\oExptCtgF_{y}(x)$ for each $x \in \mathcal{X}$.
\end{theorem}
Rather than having to calculate integrals directly to compare the performance of different policies, Theorem~\ref{theorem:policy-switch-improvement} allows us to check a uniform condition involving similarly defined functional equations.

\section{Revenue Optimization With Stochastic Prices}\label{section:multistage-revenue-optimisation}
To gain a deeper understanding of the differences between sample-based SDP and MPC, we study a particular problem from within the setting considered in Section~\ref{section:stochastic-optimal-control}. This problem is an example from the class for which sample-based MPC satisfies the DRO interpretation of Theorem~\ref{theorem:mpc-is-dro}. For initial inventory level $x^1\in\reals_+$ we consider the revenue optimization problem
\begin{alignat}{2}\label{problem:rosp-sum}
& \kern-0.07em \maximize_{y^1, y^2, \ldots} && \quad \Expt_{\Prob^{\infty}}\Biggl[\, \sum_{t=1}^{\infty}\beta^{t-1}\bigl(\rp^t(x^t-x^{t+1})-\storageCost(x^{t+1})\bigr) \Biggr] \tag{ROSP}\\
& \subjectTo &&
\quad x^{t+1} = y^{t}(x^t,\rp^t) \in [0,x^t] \quad  \forall t\in \nats. \notag
\end{alignat}
As with the general problem (\ref{problem:soc-sum}), a decision at stage $t$ is made for the next inventory level $x^{t+1}$ given the current inventory level $x^t$ and the realisation of the random variable $\rp^t$. Here $\rp^{1},\rp^{2},\ldots$ are random prices that are independent and identically distributed according to $\Prob$, and ${\storageCost:\reals_+\to\reals_+}$ is an inventory-storage cost function. We assume that $C$ is increasing, strictly convex, and continuously differentiable. Writing $\storageCostDrv(x) \defeq \frac{\drv}{\drv x} \storageCost(x)$, with $\storageCostDrv$ strictly increasing and continuous, we further define a continuous inverse function $\storageCostDrv^{-1}$ on the range of $\storageCostDrv$. Since the inventory levels are restricted to the compact set $[0,x^1]$, without loss of generality we also assume that $\lim_{x\to\infty}\storageCostDrv(x)=\infty$.

The problem (\ref{problem:rosp-sum}) is that facing a merchant who maximizes their expected discounted reward by at each time $t$ selling down from their current inventory level $x^{t}$ to a new inventory level $x^{t+1}$ at a realisation of the random price $\rp^{t}$, while incurring a storage cost ${\storageCost(x^{t+1})}$ on their remaining inventory. This can model a number of situations. For example, it could describe an electricity distributor with a charged battery deciding when to dispatch electricity, or an investor deciding when to sell a holding of shares.

As for (\ref{problem:soc-sum}), the problem (\ref{problem:rosp-sum}) is closely related to that of finding a function $\ctgF:\reals_+\times\reals\to\reals$ which satisfies the functional equation
\begin{equation}
\ctgF(x,p)=\sup_{0\leq y\leq x}\Bigl\{ p(x-y)-\storageCost(y)+\beta \Expt_{\Prob}\bigl[\ctgF(y,\rp)\bigr]\Bigr\} \quad \forall (x,p)\in{\reals_+\times\reals}.
\label{equation:rosp-recursion}
\end{equation} 
The constraints in (\ref{problem:rosp-sum}) restrict states to the compact set $[0,x^1]$, and with $\Prob$ having compact support, the terms in the sum in the objective function are bounded. It follows that (\ref{equation:rosp-recursion}) has a unique continuous solution $\ctgF:\reals_+\times\reals\to\reals$ and that the mapping $x\mapsto v(x,p)$ is concave for every $p\in\reals$ \citep[Chapter~9]{stokey-lucas:recursive-methods}. Given such $\ctgF$, a policy $y:\reals_+\times\reals\to\reals_+$ satisfying
\begin{equation}
y(x,p) \in \argmax_{0\leq y\leq x}\Bigl\{ p(x-y)-\storageCost(y)+\beta \Expt_{\Prob}\bigl[\ctgF(y,\rp)\bigr]\Bigr\} \quad \forall (x,p)\in{\reals_+\times\reals}
\label{equation:rosp-argmax}
\end{equation}
is optimal for (\ref{problem:rosp-sum}). 

In the notation of Section~\ref{section:stochastic-optimal-control}, (\ref{problem:rosp-sum}) has $\mathcal{Y}(x,p) = [0,x]$, and a suitable change of perspective from maximization to minimization yields the cost function $\varphi(x,y,p) = -p(x-y)+\storageCost(y)$ which is concave in $p$. This satisfies the conditions of Proposition~\ref{proposition:concave-costs}, so applying sample-based MPC to (\ref{problem:rosp-sum}) will yield the DRO interpretation of Theorem~\ref{theorem:mpc-is-dro}.

Denote the projection of $z \in \reals$ onto $[a,b]\subset\reals$ by $(z)_{[a,b]}$ and write $(z)_{+}$ for its projection onto the nonnegative real numbers. Using the inclusion (\ref{equation:rosp-argmax}) we arrive at the following proposition which presents the optimal policy in closed form.
\begin{proposition}
Suppose that the distribution $\Prob$ has compact support. Then the policy which at inventory $x$ and price $p$ sells down to
\begin{equation*}
y(x,p) = \storageCostDrv^{-1}\Bigl(\bigl(\beta\Expt_{\Prob}\bigl[(\rp-{p})_+\bigr]-(1-\beta)p\bigr)_{[\storageCostDrv(0),\storageCostDrv(x)]} \Bigr)
\end{equation*}
is optimal for (\ref{problem:rosp-sum}).
\label{proposition:rosp-optimal-policy}
\end{proposition}
The optimal policy has a natural interpretation: the term $\beta  \Expt_{\Prob}\bigl[(\rp-{p})_+\bigr]$ is the discounted expected increase in price gained by not selling, and the term $(1-\beta)p$ is the portion of the current price that is lost due to discounting by not selling. The difference in these terms gives the expected net increase in price gained from storing inventory, and this is balanced against the marginal storage cost.

Proposition~\ref{proposition:rosp-optimal-policy} shows that for each price $p$ the problem (\ref{problem:rosp-sum}) has an optimal target inventory level $y(\infty,p) \defeq \lim_{x\to\infty}y(x,p)$. At inventory $x$ the optimal policy is to sell down to $y(\infty,p)$ if $x$ is above this, and sell nothing otherwise. Conversely, there is a minimum acceptable price $\cPrice(x)$ required for selling any portion of inventory from $x$ to be worthwhile: this is the highest price $p$ which solves $y(\infty,p)=x$. Note that due to the continuity of $c^{-1}$, the function $\cPrice$ is also continuous.

Proposition~\ref{proposition:rosp-optimal-policy} makes no assumptions about the distribution $\Prob$, except that it has compact support. Thus, it could have a density on a compact set, or could be an empirical distribution. Using the empirical distribution on $N$ price samples ${p}_{1},\ldots,{p}_{N}$, Proposition~\ref{proposition:rosp-optimal-policy} shows that the optimal SDP policy, which we denote by $y_{\SDP}$, is
\begin{equation}
y_{\SDP}(x,p;p_1,\ldots,p_N)=\storageCostDrv^{-1}\Bigl(\bigl(\beta\tfrac{1}{N}\tsum_{i=1}^{N}({p}_{i}-{p})_{+}-(1-\beta)p\bigr)_{[\storageCostDrv(0),\storageCostDrv(x)]}\Bigr).
\label{equation:sdp-policy}
\end{equation}
The optimal MPC policy, which we denote by $y_{\MPC}$, can be then be obtained from Proposition~\ref{proposition:rosp-optimal-policy} by applying it to the point-mass distribution at the sample average $\mu_N = \sum_{i=1}^{N}{p}_{i}$, giving
\begin{equation}
y_{\MPC}(x,p;p_1,\ldots,p_N)=\storageCostDrv^{-1}\Bigl(\bigl(\beta(\mu_N-{p})_{+}-(1-\beta){p}\bigr)_{[\storageCostDrv(0),\storageCostDrv(x)]}\Bigr).  
\label{equation:mpc-policy}
\end{equation}
In what follows, for simplicity we often suppress the dependence on the price samples $p_1,\ldots,p_N$. We denote by $y_{\SDP}(\infty,p)$ and $y_{\MPC}(\infty,p)$ the target inventory levels of the policies (\ref{equation:sdp-policy}) and (\ref{equation:mpc-policy}), and by $\cPriceSDP(x)$ and $\cPriceMPC(x)$ their minimum acceptable prices required for sales. Clearly it holds that $\cPriceSDP \leq \max\{p_1,\ldots,p_N\}$, since this is the highest price the SDP policy considers a possibility. Similar remarks hold for the MPC policy.

Depending on the values of the price samples $p_1,\ldots,p_N$, the sample-based policies may hold on to inventory for too long and overpay for storage, or sell inventory prematurely and not be able to capitalise on higher prices offered in the future. Jensen's Inequality shows that ${\sum_{i=1}^{N}({p}_{i}-{p})_{+} \geq (\mu_N-{p})_{+}}$, and thus $y_{\SDP}(\infty,p) \geq y_{\MPC}(\infty,p)$. This reveals that the MPC policy decides that it is worthwhile to sell at lower prices than the SDP policy does. Indeed, $\cPriceSDP(x) \geq \cPriceMPC(x)$ for every inventory level $x$, and the trade-off between overpaying for storage and selling inventory prematurely is handled differently by sample-based SDP and MPC.

\subsection{Out-of-Sample Performance}
\label{section:sic-out-of-sample} 
To study the out-of-sample performance of SDP and MPC on (\ref{problem:rosp-sum}), we make the following assumption regarding the true distribution.
\begin{assumption}
\label{assumption:rosp-unbounded-prices} The 
distribution $\Prob$ has nonnegative support, a finite mean, and no atoms.
\end{assumption}
Under Assumption~\ref{assumption:rosp-unbounded-prices}, for every $x,y\in[0,x^1]$ and $\Prob$-almost every $p\in\reals_+$, it holds that 
\begin{equation*}
    \bigl\lvert p(x-y)-\storageCost(y)\bigr\rvert \leq px^1+\storageCost(x^1).
\end{equation*}
With $\Expt_{\Prob}\bigl[\rp x^1+\storageCost(x^1)\bigr]$ finite valued, this shows that Assumption~\ref{assumption:rosp-unbounded-prices} for (\ref{problem:rosp-sum}) is akin to Assumption~\ref{assumption:unbounded-noise} for (\ref{problem:soc-sum}), which enables us to define the out-of-sample performance of policies using functional equations. 

Under Assumption~\ref{assumption:rosp-unbounded-prices} the out-of-sample performance of the SDP policy (\ref{equation:sdp-policy}) when starting from initial inventory $x$, which we denote by $\oExptCtgF_{\SDP}(x)$, is well defined and satisfies the functional equation 
\begin{equation}\label{equation:rosp-oos-sdp-value-function}
\oExptCtgF_{\SDP}(x)=\Expt_{\Prob}\Bigl[ \rp\bigl(x-y_{\SDP}(x,\rp)\bigr)-\storageCost\bigl(y_{\SDP}(x,\rp)\bigr)+\beta \oExptCtgF_{\SDP}\bigl(y_{\SDP}(x,\rp)\bigr)\Bigr].
\end{equation}
Similarly, the out-of-sample performance of the MPC policy (\ref{equation:mpc-policy}) when starting from initial inventory $x$, which we denote by $\oExptCtgF_{\MPC}(x)$, is well defined and satisfies the functional equation 
\begin{equation}\label{equation:rosp-oos-mpc-value-function}
\oExptCtgF_{\MPC}(x)=\Expt_{\Prob}\Bigl[\rp\bigl(x-y_{\MPC}(x,\rp)\bigr)-\storageCost\bigl(y_{\MPC}(x,\rp)\bigr)+\beta \oExptCtgF_{\MPC}\bigl(y_{\MPC}(x,\rp)\bigr)\Bigr].
\end{equation}

Let $\Indicator$ be the event indicator; that is $\Indicator\{ \mathcal{E}\}=1$ if the event $\mathcal{E}$ is true and $\Indicator \{\mathcal{E}\}=0$ if it is false. We use Theorem~\ref{theorem:policy-switch-improvement} to establish the following proposition which compares the out-of-sample performance of SDP and MPC.
\begin{proposition}\label{proposition:MPC-wins-with-outliers}
Let Assumption~\ref{assumption:rosp-unbounded-prices} hold. For price samples $p_1,\ldots,p_N$ which determine the sample-based SDP and MPC policies by (\ref{equation:sdp-policy}) and (\ref{equation:mpc-policy}), if the minimum acceptable SDP price $\cPriceSDP(x;p_1,\ldots,p_N)$ is such that
\begin{equation*}
    \storageCostDrv(x) \geq \beta \Expt_{\Prob}\Bigl[\,\rp\Indicator\bigl\{\rp > \cPriceSDP(x;p_1,\ldots,p_N)\bigr\}\Bigr]  \quad \forall x\in [0,x^1], 
\end{equation*}
then $\oExptCtgF_{\MPC}(x;p_1,\ldots,p_N) \geq \oExptCtgF_{\SDP}(x;p_1,\ldots,p_N)$ for each $x\in [0,x^1]$. That is, the MPC policy performs at least as well out of sample on (\ref{problem:rosp-sum}) as the SDP policy.
\end{proposition}
The condition has a natural interpretation: $\beta\Expt_{\Prob}\bigl[\rp\Indicator\{\rp > \cPriceSDP(x;p_1,\ldots,p_N)\}\bigr]$ represents the expected revenue that SDP gains over MPC by delaying sales, and if this is less than the marginal cost of storage, then waiting provides no net benefit. Without loss of generality, assume ${p}_1 \leq \cdots \leq {p}_N$ and recall that $\cPriceSDP (x;p_1,\ldots,p_N)\leq {p}_N$. The mapping $p_N \mapsto ({p}_N-{p})_{+}$ in (\ref{equation:sdp-policy}) is strictly increasing at the point ${p}=\cPriceSDP(x;p_1,\ldots,p_N)$, where $y_{\SDP}(\infty,p;p_1,\ldots,p_N)=x$. It follows that $p_N \mapsto \cPriceSDP(x;p_1,\ldots,p_N)$ is strictly increasing, and thus $p_N \mapsto \Expt_{\Prob}\bigl[\rp\Indicator\{\rp > \cPriceSDP(x;p_1,\ldots,p_N)\}\bigr]$ is decreasing. If the $N$ price samples are independent and identically distributed according to $\Prob$, when $\Prob$ has a small amount of probability at very high prices, we will occasionally sample a maximal price that is sufficiently large for the inequality to hold.

Proposition~\ref{proposition:MPC-wins-with-outliers} is not an explicit statement about the relative expected out-of-sample performance of SDP and MPC under the sampling distribution. In fact, in most applications it is likely that there will always be some outcomes of the samples which result in MPC outperforming SDP, so Proposition~\ref{proposition:MPC-wins-with-outliers} is not surprising. However, the result demonstrates how the relative performance of SDP and MPC is affected by the form of true price distribution: the sample-based SDP policy can be misled by overly high price samples that cause it to hold on to inventory for too long and overpay for storage. Consequently, right skew and the size of the right-hand tail in the true distribution are likely to impact the performance of SDP. On the other hand, due to its distributionally robust properties, the sample-based MPC policy more quickly sells down to lower inventory levels and is protected against this.

If $\storageCostDrv(0)=0$, Proposition~\ref{proposition:MPC-wins-with-outliers} does not apply. But this is to be expected; when $\storageCostDrv(0)=0$ an infinitesimal amount of inventory incurs negligible storage costs, so waiting longer under the sample-based SDP policy to observe a higher price will always perform better in the long run. Despite this, for nonnegligible initial inventory levels, we present examples below which show that MPC can still outperform SDP when $\storageCostDrv(0)=0$, due to high storage costs accumulated early on.

\subsection{Expected Out-of-Sample Performance}
In the rest of this section we present examples which compare the {expected} out-of-sample performance of SDP and MPC on (\ref{problem:rosp-sum}) for price distributions with different skews and tail sizes. We suppose that the $N$ price samples $\rp_1,\ldots,\rp_N$ used to determine the sample-based policies by (\ref{equation:sdp-policy}) and (\ref{equation:mpc-policy}) are independent and identically distributed according to $\Prob$. Expected out-of-sample performance is then given by the expectations of (\ref{equation:rosp-oos-sdp-value-function}) and (\ref{equation:rosp-oos-mpc-value-function}) under the sampling distribution of $\rp_1,\ldots,\rp_N$.

Let $\text{Exponential}(\lambda)$ denote the exponential distribution with rate $\lambda$. Exponential distributions are strongly right skewed. When prices are exponentially distributed, we have the following result.
\begin{proposition}\label{proposition:exponential-example} 
    Let $\Prob$ be the $\text{\normalfont Exponential}(1)$ distribution, $x^1 = 1$, and $\storageCost(x) = \frac{1}{2}x^{2}$. For each sample size $N \geq 2$, as the discount factor $\beta\to 1$ the expected out-of-sample performance of SDP on (\ref{problem:rosp-sum}) is unbounded below, while that of MPC is bounded. 
\end{proposition}

Proposition~\ref{proposition:exponential-example} shows that for each sample size greater than or equal to $2$ the expected out-of-sample performance of SDP can be made arbitrarily worse than that of MPC by choosing a discount factor that is sufficiently close to $1$. For a particular discount factor the performance of SDP could be improved by increasing the sample size, but this may be required to be very large for SDP to outperform MPC. The result shows that the performance of SDP can be arbitrarily bad even for arbitrarily large sample sizes. Moreover, in real applications the sample-based SDP policy must be computed numerically, but algorithms for doing this cannot handle large sample sizes. MPC does not have this problem.

\bigskip

To compute the expected out-of-sample performance of SDP and MPC on (\ref{problem:rosp-sum}) for different price distributions, we use a simulation coded in \texttt{Julia} \citep{bezanson2017julia}. Although (\ref{problem:rosp-sum}) has an infinite horizon, simulation with a finite number of stages (say $T$) will provide accurate performance estimates as long as this is sufficiently large; we set $T = 1\cdot 10^3$.\footnote{The resulting truncation error has an order less than $0.99^{1000}/(1-0.99) \approx 1\cdot 10^{-3}$.} With $x^{1} = 1$, $\beta = 0.99$, and $\storageCost(x) = \frac{1}{2}x^{2}$, for a given sample size $N$ we repeat:
\begin{enumerate}[label={\normalfont\arabic*.}, ref={\arabic*.}]
\setcounter{enumi}{-1}
\item\label{algorithm:compute-out-of-sample-cost-1} Sample $\rp_1,\ldots,\rp_N\sim\Prob$ to determine a policy $y$ using (\ref{equation:sdp-policy}) for SDP and (\ref{equation:mpc-policy}) for MPC.
\item\label{algorithm:compute-out-of-sample-cost-2} Sample $\rp^{t}\sim\Prob$, compute ${\beta^{t-1}\bigl(\rp^{t}(x^{t}-y(x^{t},\rp^{t}))-\storageCost(y(x^{t},\rp^{t}))\bigr)}$, and set $x^{t+1} = y(x^{t},\rp^{t})$.
\item\label{algorithm:compute-out-of-sample-cost-3} Repeat \ref{algorithm:compute-out-of-sample-cost-2} from stage $t=1$ to $T$ and compute ${\sum_{t=1}^{T} \beta^{t-1}\bigl(\rp^{t}(x^t-y(x^{t},\rp^{t}))-\storageCost(y(x^{t},\rp^{t}))\bigr)}$.
\end{enumerate}
\noindent Each repetition of \ref{algorithm:compute-out-of-sample-cost-1}--\ref{algorithm:compute-out-of-sample-cost-3} generates a random out-of-sample performance realisation, and the average of these realisations provides a statistical estimate of expected out-of-sample performance. In the following simulations we used $1\cdot 10^5$ realisations, finding this sufficient to ensure accurate results. (In Figures~\ref{figure:left-triangle-sim}--\ref{figure:lognormal-sim} the standard error ranges are smaller than the markers so are not shown.)

\subsection{Skewed Price Distributions}
\begin{figure}[H]
        \centering
        \hspace{-1.5cm}\includegraphics[width=289pt]{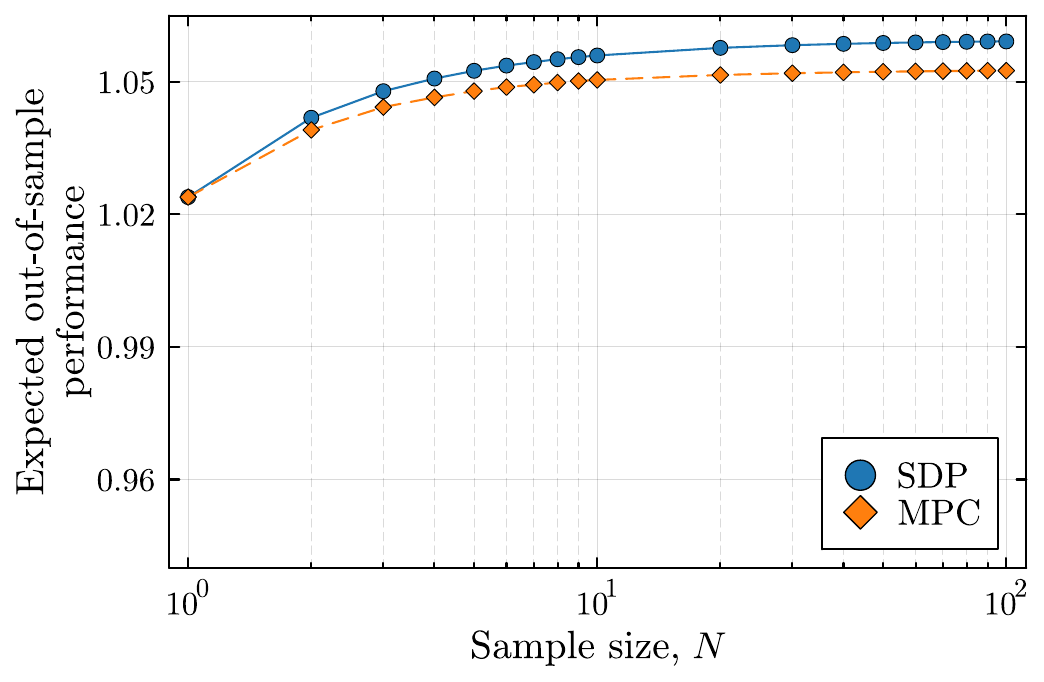}
        \caption{Performance of SDP and MPC on (\ref{problem:rosp-sum}) for $\bm{\text{Triangular}(0,}\nicefrac{\bm{3}}{\bm{2}}\bm{,}\nicefrac{\bm{3}}{\bm{2}}\bm{)}$ distributed price.}
        \label{figure:left-triangle-sim}
\end{figure}

\begin{figure}[H]
        \centering
        \hspace{-1.5cm}\includegraphics[width=289pt]{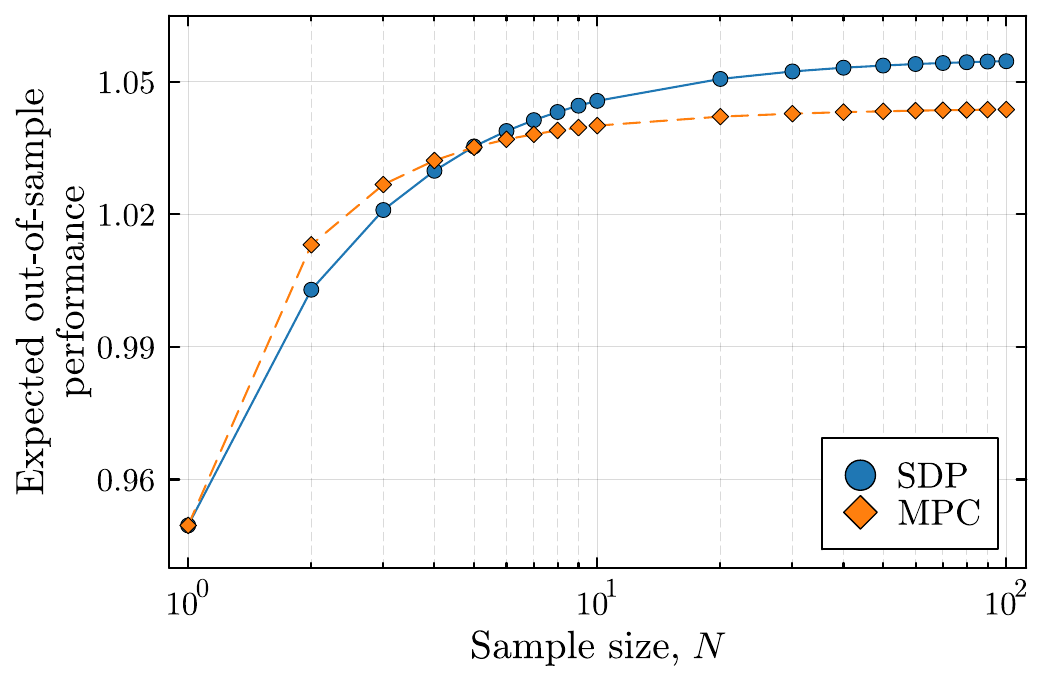}
        \caption{Performance of SDP and MPC on (\ref{problem:rosp-sum}) for $\bm{\text{Triangular}(}\nicefrac{\bm{1}}{\bm{2}}\bm{,}\nicefrac{\bm{1}}{\bm{2}}\bm{,2)}$ distributed price.}
        \label{figure:right-triangle-sim}
\end{figure}

Let $\text{Triangular}(a,m,b)$ denote the triangular distribution with lower limit $a$, mode $m$, and upper limit $b$. Figures~\ref{figure:left-triangle-sim} and \ref{figure:right-triangle-sim} present the expected out-of-sample performance of SDP and MPC for left- and right-skewed triangular price distributions with mean $1$ and variance $\nicefrac{1}{8}$.

\noindent Figure~\ref{figure:left-triangle-sim} shows SDP outperforming MPC for all sample sizes when the true price distribution is left skewed. In contrast, Figure~\ref{figure:right-triangle-sim} shows MPC outperforming SDP for sample sizes $N \leq 4$ when the true price distribution is right skewed. This reflects the fact that the likelihood of the samples containing a price high enough to cause SDP to undersell inventory and overpay for storage is more likely in the presence of right skew.

\begin{figure}[H]
        \centering
        \hspace{-1.5cm}\includegraphics[width=289pt]{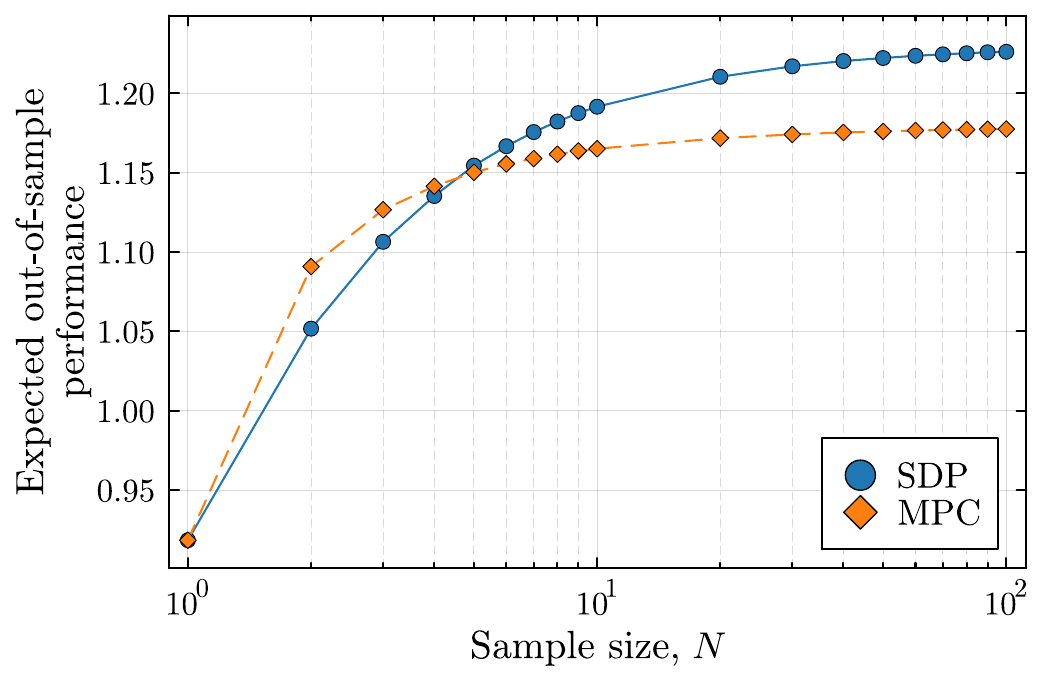}
        \caption{Performance of SDP and MPC on (\ref{problem:rosp-sum}) for $\bm{\text{Triangular}(0,0,3)}$ distributed price.}
        \label{figure:right-triangle-higher-variance-sim}
\end{figure}

Figure~\ref{figure:right-triangle-higher-variance-sim} presents the expected out-of-sample performance of SDP and MPC for a right-skewed triangular price distribution with mean $1$ and variance $\nicefrac{1}{2}$.
\noindent Figure~\ref{figure:right-triangle-higher-variance-sim} shows MPC outperforming SDP for sample sizes $N \leq 4$ when the true price distribution is right skewed and has a higher variance than that used in Figure~\ref{figure:right-triangle-sim}. This is the same range as that in Figure~\ref{figure:right-triangle-sim}, but here the relative amount by which MPC outperforms SDP is increased.

\subsection{Right-Tailed Price Distributions}
\begin{figure}[H]
        \centering
        \hspace{-1.5cm}\includegraphics[width=289pt]{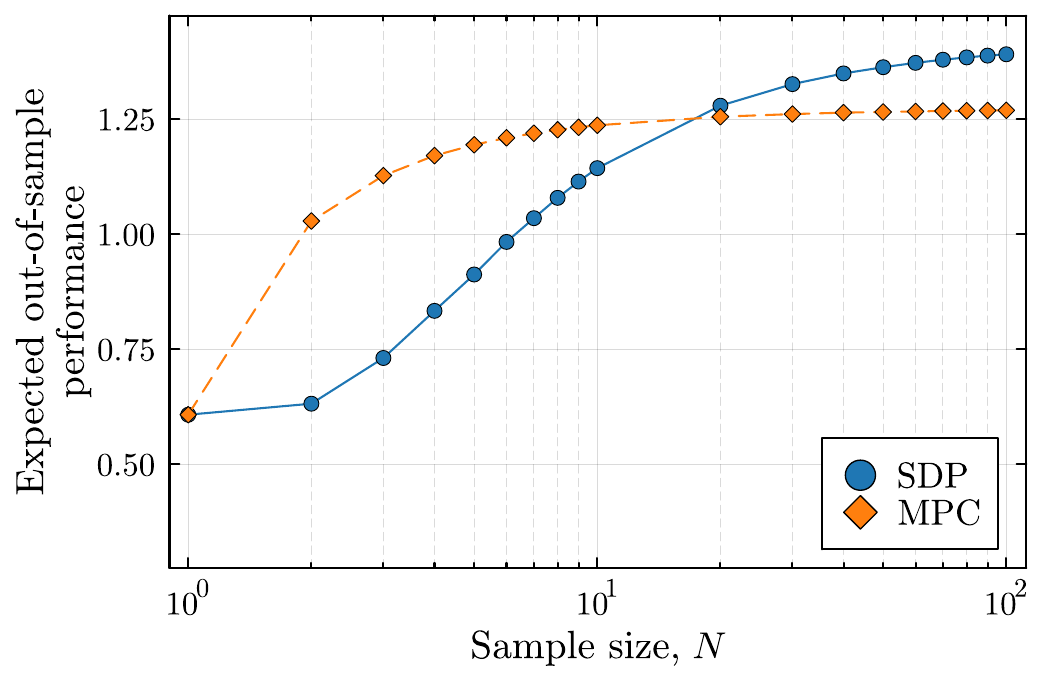}
        \caption{Performance of SDP and MPC on (\ref{problem:rosp-sum}) for $\bm{\text{Exponential}(1)}$ distributed price.}
        \label{figure:exponential-sim}
\end{figure}

\begin{figure}[H]
        \centering
        \hspace{-1.5cm}\includegraphics[width=289pt]{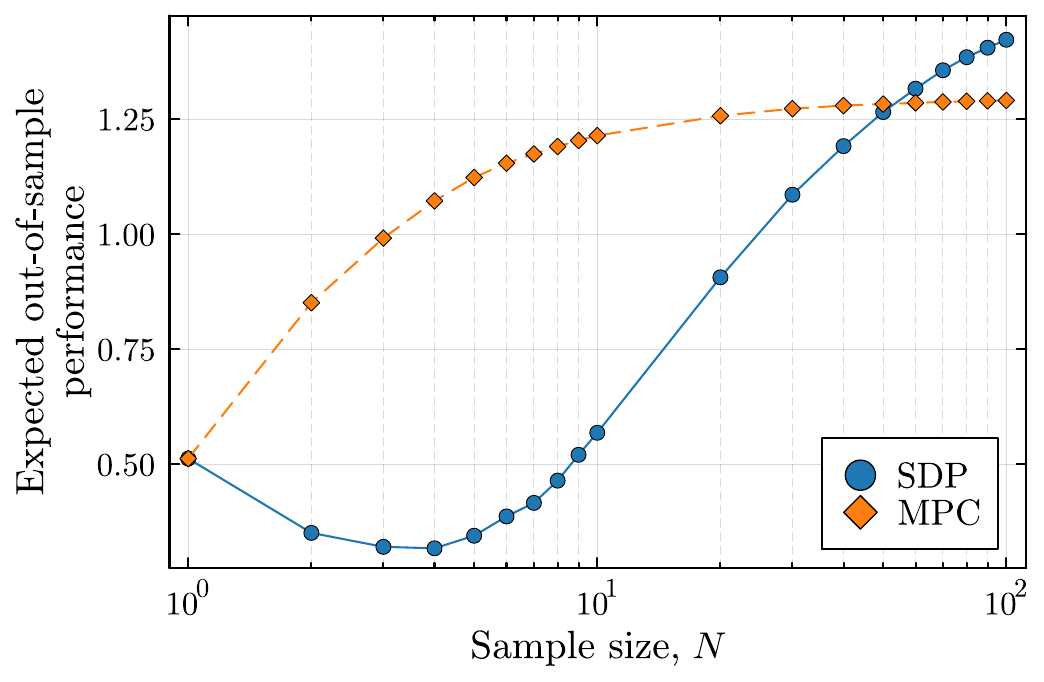}
        \caption{Performance of SDP and MPC on (\ref{problem:rosp-sum}) for $\bm{\text{LogNormal}(-}\nicefrac{\bm{1}}{\bm{2}}\bm{,1)}$ distributed price.}
        \label{figure:lognormal-sim}
\end{figure}

Recall that $\text{Exponential}(\lambda)$ is the exponential distribution with rate $\lambda$, and let $\text{LogNormal}(\mu,\sigma^2)$ denote the lognormal distribution with mean  $\mu$ and variance $\sigma^2$. Figures~\ref{figure:exponential-sim} and \ref{figure:lognormal-sim} present the expected out-of-sample performance of SDP and MPC for exponential and lognormal price distributions with mean $1$.
\noindent Figure~\ref{figure:exponential-sim} shows MPC outperforming SDP for sample sizes $N \leq 10$ when the true price distribution is exponential. This is a larger range than those in Figures~\ref{figure:left-triangle-sim}--\ref{figure:right-triangle-higher-variance-sim}. Here the right-hand tail of the exponential distribution increases the propensity for overly high prices to be included in the samples, which worsens the relative performance of SDP. 
\noindent Figure~\ref{figure:lognormal-sim} shows MPC outperforming SDP for samples sizes $N \leq 50$ when the true price distribution is right skewed and has a heavier right-hand tail the exponential distribution used in Figure~\ref{figure:exponential-sim}. The poor performance of SDP exhibited here suggests that a result similar to that of Proposition~\ref{proposition:exponential-example} for exponential distributions may also hold for lognormal distributions. Considering that lognormal distributions have more weight in their tails than exponential distributions, this would not be surprising.

\begin{figure}[H]
\centering
\hspace{0.9cm}\includegraphics[width=281pt]{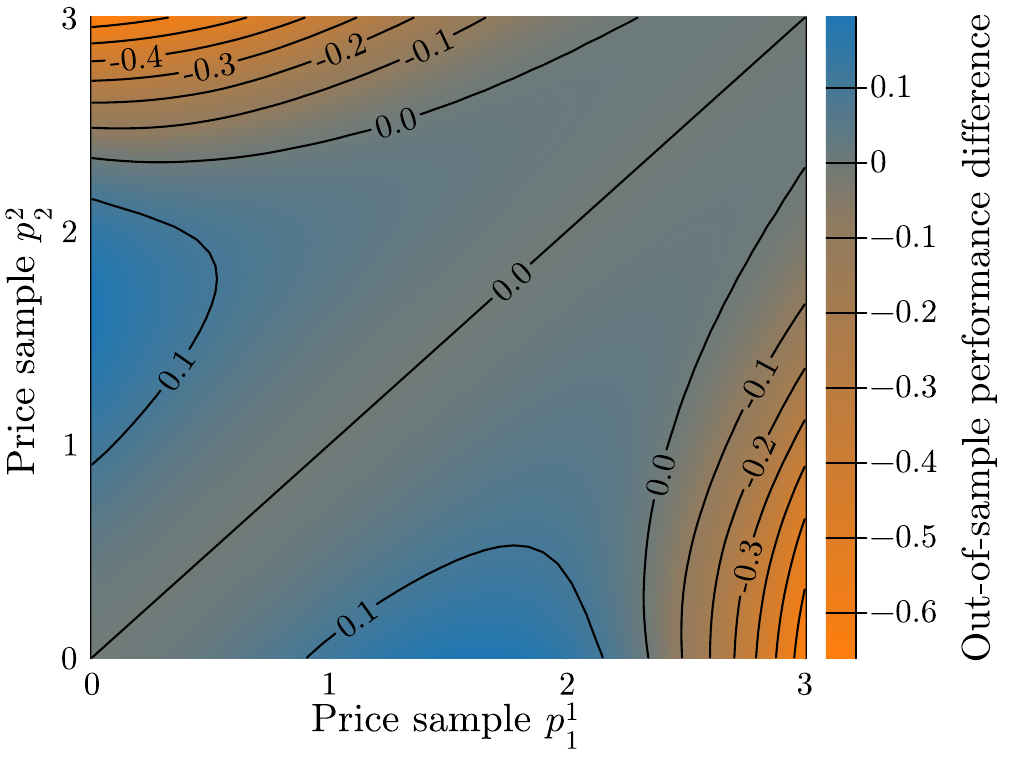}
\caption{\textbf{Difference in SDP and MPC performance on (\ref{problem:rosp-sum}) for $\bm{\text{LogNormal}(-}\nicefrac{\bm{1}}{\bm{2}}\bm{,1)}$ distributed price and $\bm{N=2}$ samples. Positive values indicate when SDP outperforms MPC and negative values indicate the reverse.}}
\label{figure:log-normal-contours}
\end{figure}

For a sample of size $N=2$, Figure~\ref{figure:log-normal-contours} presents the difference in the out-of-sample performance of the SDP and MPC policies.

\vspace{-2.77cm}\hspace{8.83cm}\textcolor{white}{\Large{$\bullet$}}\vspace{2.77cm}

\vspace{-8.32cm}\hspace{3.37cm}\textcolor{white}{\Large{$\bullet$}}\vspace{8.32cm}

\vspace{-1.3cm}

\noindent Figure~\ref{figure:log-normal-contours} explicitly demonstrates that a single overly high price sample worsens the performance of SDP more than that of MPC; typical samples result in the SDP policy outperforming the MPC policy, but for extremes where one sample is very high, the reverse occurs, as predicted by Proposition~\ref{proposition:MPC-wins-with-outliers}.

\section{Discussion}
\label{section:SDP-vs-MPC-discussion}
In this paper we provide an explanation for the good out-of-sample performance of MPC that is sometimes observed in practice, based on the interpretation of MPC as solving a mean-constrained distributionally robust optimization problem. This depends critically on the MPC Bellman operator associated with the problem being concavity preserving. When the operator is instead convexity preserving, the interpretation of MPC becomes distributionally optimistic. Loosely speaking, the problem type for which the robust interpretation holds has uncertainty in the objective function parameters \citep[as in][]{downward_et_al_2020_objective_uncertainty}, and that for which the optimistic interpretation holds features constraint uncertainty \citep[as in][]{pereira1991multi}. These features provide some guidance on when MPC may or may not be the right choice to address uncertainty within a particular problem. However, this is not the whole story, since out-of-sample performance also depends on the form of the underlying probability distributions present in a particular problem.

We study the out-of-sample performance of SDP and MPC on a stochastic revenue optimization problem of the type in \citep{downward_et_al_2020_objective_uncertainty}, finding that SDP can be outperformed by MPC when the underlying price distribution is right skewed or has a large right-hand tail and the number of samples is not too large. In the case where the underlying price distribution is exponential and the discount factor approaches $1$, SDP can be outperformed by MPC regardless of the number of samples. 

The issues which occur when applying sample-based SDP to the stochastic revenue optimization problem may be alleviated by appending newly observed prices to the sample history and updating the policy before applying it again. In general, this approach is not practical as the time complexities of algorithms used to solve for the optimal (in-sample) SDP policy grow quickly in the number of samples. It could be more effective to use some sort of rolling window so that policies are based on a subset of the most recent samples. Still, this requires recalculation of the optimal policy at each stage, which may not be desirable. Also notice that whenever such a rolling window includes a very high price sample, the same performance issues will occur.

The stochastic revenue optimization problem we consider is quite restricted, having deterministic dynamics and constraints that do not allow inventory to increase. This gives a transient problem of selling inventory down rather than a problem with a steady-state component. We have chosen this problem because it enables an analytical solution and a more detailed analysis, but we expect that similar results would be obtained for problems in which there are occasional additional amounts of inventory arriving. For example, if stock to replenish inventory to an initial level arrives in each time period with probability $\gamma$, then a renewal-theory argument shows that the problem of maximizing average reward per unit time is equivalent to the revenue optimization problem with the discount factor $\beta = 1- \gamma$.

%
%
%
\begin{APPENDICES}

\section{Proof of Proposition~\ref{proposition:concave-costs}}
\begin{repeatproposition}[Proposition~\ref{proposition:concave-costs}.]
Let Assumption~\ref{assumption:bounded-continuous} hold. Suppose that the set-valued mapping $\mathcal{Y}$ does not depend on $\xi$ and that the mapping $\xi\mapsto\stageC(x,y,\xi)$ is concave for every $x, y \in \mathcal{X}$. Then the MPC Bellman operator $B_{\MPC}$ is concavity preserving.
\end{repeatproposition}

\proof{Proof of Proposition~\ref{proposition:concave-costs}}
    Consider the set
\begin{equation*}
    \mathcal{F} = \bigl\{f\in\bcfcns(\mathcal{X}\times\Xi) : \xi\mapsto f(x,\xi)~\text{is concave for every}~x\in\mathcal{X}\bigr\}.
\end{equation*}
By definition, $\xi\mapsto f(x,\xi)$ is concave for every $x\in\mathcal{X}$ and $f\in\mathcal{F}$. It remains to show that for each $f\in\mathcal{F}$, we have $B_{\MPC}(f)\in\mathcal{F}$. Let $f\in\mathcal{F}$ and observe that for every $x\in\mathcal{X}$, $\xi, \xi^{\prime}\in\Xi$, and $\lambda\in[0,1]$, it holds that 
\begin{align*}
    B_{\MPC}(f)\bigl(x,\lambda\xi+(1-\lambda)\xi^{\prime}\bigr) &= \inf_{y\in\mathcal{Y}(x)}\Bigl\{\stageC\bigl(x,y,\lambda\xi+(1-\lambda)\xi^{\prime}\bigr)+\beta f(y,\mu_N )\Bigr\}\\
    &= \stageC\bigl(x,y^{\star},\lambda\xi+(1-\lambda)\xi^{\prime}\bigr)+\beta f(y^{\star},\mu_N )\\
    &\geq \lambda\bigl(\stageC(x,y^{\star},\xi)+\beta f(y^{\star},\mu_N )\bigr) + (1-\lambda)\bigl(\stageC(x,y^{\star},\xi^{\prime})+\beta f(y^{\star},\mu_N )\bigr)\\
    &\geq \lambda \inf_{y\in\mathcal{Y}(x)}\Big\{\stageC(x,y,\xi)+\beta f(y,\mu_N )\Big\} + (1-\lambda)\inf_{y\in\mathcal{Y}(x)}\Big\{\stageC(x,y,\xi^{\prime})+\beta f(y,\mu_N )\Big\}\\
    &=\lambda B_{\MPC}(f)(x,\xi) + (1-\lambda)B_{\MPC}(f)(x,\xi^{\prime}), 
\end{align*} 
where the second equality holds for some $y^{\star}\in\mathcal{Y}(x)$ since the infimum of a continuous function over a compact set is attained, the first inequality holds due to the concavity of the mapping $\xi\mapsto\varphi(x,y^{\star},\xi)$, and the last inequality holds since $y^{\star}\in\mathcal{Y}(x)$. This implies that $\xi\mapsto B_{\MPC}(f)(x,\xi)$ is concave for every $x\in\mathcal{X}$ and hence that $B_{\MPC}(f)\in\mathcal{F}$.\Halmos
\endproof
\medskip

\section{Proof of Theorem~\ref{theorem:mpc-is-dro} and Corollary~\ref{corollary:mpc-is-doo}}

\begin{repeattheorem}[Theorem~\ref{theorem:mpc-is-dro}.]
Let Assumption~\ref{assumption:bounded-continuous} hold. For an ambiguity set $\mathcal{P}$, suppose that $\Expt_{\Qrob}[\rxi] = \mu_N$ for every $\Qrob\in\mathcal{P}$ and that $\pmProb_{\mu_N} \in \mathcal{P}$. If the MPC Bellman operator $B_{\MPC}$ is concavity preserving, then the DRO functional equation (\ref{equation:dro-value-function}) has a solution which is the same as the solution to the MPC functional equation (\ref{equation:mpc-value-function}).
\end{repeattheorem}

\proof{Proof of Theorem~\ref{theorem:mpc-is-dro}}
The equation (\ref{equation:dro-value-function}) features a DRO Bellman operator $B_{\DRO}$ which for ${f \in \bcfcns(\mathcal{X}\times\Xi)}$ at $(x,\xi)\in\mathcal{X}\times\Xi$ has value
\begin{equation*}
    B_{\DRO}(f)(x,\xi) = \inf_{y\in\mathcal{Y}(x,{\xi})}\Bigl\{\stageC(x,y,{\xi})+\beta \sup_{\Qrob\in\mathcal{P}}\Expt_{\Qrob}\bigl[f(y,\rxi)\bigr]\Bigr\}.
\label{equation:dro-bellman-operator}
\end{equation*}
With $B_{\MPC}$ concavity preserving, there exists a closed subset $\mathcal{F}\subseteq\bcfcns(\mathcal{X}\times\Xi)$ where the mapping $\xi\mapsto f(x,\xi)$ is concave for every $x\in\mathcal{X}$ and $f\in \mathcal{F}$, and which maps to itself under $B_{\MPC}$. We first show that $B_{\DRO}$ is a contraction mapping on $\mathcal{F}$ via comparison to $B_{\MPC}$.

Let $y\in\mathcal{X}$. For $f\in\mathcal{F}$ we have that
\begin{equation}
\label{equation:sup-upper-bound}
     \Expt_{\Qrob}\bigl[f(y,\rxi)\bigr] \leq f(y,\Expt_{\Qrob}[\rxi]) = f(y,\mu_N ) \quad \forall \Qrob\in\mathcal{P}, 
\end{equation}
where the inequality follows from Jensen's inequality (due to concavity), and the equality follows since ${\Expt_{\Qrob}[\rxi] = \mu_N}$ for all $\Qrob \in \mathcal{P}$. The statement (\ref{equation:sup-upper-bound}) implies that $\sup_{\Qrob \in \mathcal{P}} \Expt_{\Qrob}\bigl[f(y,\rxi)\bigr] \leq f(y,\mu_N )$. But, with $\Expt_{\pmProb_{\mu_N}}\bigl[f(y,\rxi)\bigr] = f(y,\mu_N)$ and $\pmProb_{\mu_N}\in\mathcal{P}$, this supremum is attained by $\pmProb_{\mu_N}$ and $\sup_{\Qrob \in \mathcal{P}} \Expt_{\Qrob}\bigl[f(y,\rxi)\bigr] = f(y,\mu_N )$. Thus, for every $f\in\mathcal{F}$, we have that
\begin{align}\label{equation:dro-mpc-agree-pointwise}
    B_{\DRO}(f)(x,\xi) &= \inf_{y\in\mathcal{Y}(x,\xi)}\Bigl\{\stageC(x,y,\xi)+\beta 
    \sup_{Q\in\mathcal{P}} \Expt_{\Qrob} \bigl[f(y,\rxi)\bigr]\Bigr\}\notag\\
    &= \inf_{y\in\mathcal{Y}(x,\xi)}\Bigl\{\stageC(x,y,\xi)+\beta f(y,\mu_N )\Bigr\} = B_{\MPC}(f)(x,\xi) \quad \forall (x,\xi)\in\mathcal{X}\times\Xi,
\end{align}
i.e., $B_{\DRO}$ agrees with $B_{\MPC}$ on $\mathcal{F}$. Recalling that $B_{\MPC}$ is a contraction mapping on $\mathcal{F}$ under the sup norm, the same holds for $B_{\DRO}$.

It remains to show that the DRO functional equation (\ref{equation:dro-value-function}) has a solution which is the same as the solution to the MPC functional equation (\ref{equation:mpc-value-function}). Observe that $(\mathcal{F},\lVert\blank\rVert_{\infty})$ is a Banach space since $\mathcal{F}$ is closed. As $B_{\DRO}$ is a contraction mapping on $\mathcal{F}$, Banach's fixed-point theorem shows that the equation $f = B_{\DRO}(f)$ has a unique solution in $\mathcal{F}$, and this solution solves the DRO functional equation (\ref{equation:dro-value-function}). In light of (\ref{equation:dro-mpc-agree-pointwise}), this solution is the same as the solution to (\ref{equation:mpc-value-function}).
\Halmos
\endproof

\begin{repeatcorollary}[Corollary~\ref{corollary:mpc-is-doo}.]
Suppose that the conditions in the statement of Theorem~\ref{theorem:mpc-is-dro} hold, but that the MPC Bellman operator $B_{\MPC}$ is convexity preserving  rather than concavity preserving. Then the DOO functional equation (\ref{equation:doo-value-function}) has a solution which is the same as the solution to the MPC functional equation (\ref{equation:mpc-value-function}).
\end{repeatcorollary}

\proof{Proof of Corollary~\ref{corollary:mpc-is-doo}}
The proof repeats the proof of Theorem~\ref{theorem:mpc-is-dro}, but with the inequalities reversed and $\inf_{\Qrob \in \mathcal{P}} \Expt_{\Qrob}\bigl[f(y,\rxi)\bigr]$ attained by $\pmProb_{\mu_N}$. \Halmos
\endproof
\medskip

\section{Proof of Theorem~\ref{theorem:mpc-dro-guarantee} and Corollary~\ref{corollary:mpc-doo-dissapointment}}
\begin{repeattheorem}
Let Assumption~\ref{assumption:bounded-continuous} hold, and set $\mu_N=\mu \defeq \Expt_{\Prob}[\rxi]$ in the MPC functional equation (\ref{equation:mpc-value-function}). If the resulting MPC Bellman operator $B_\MPC$ is concavity preserving, then
\begin{equation*}
    \Expt_{\Prob^\infty}\Biggl[\,\sum_{t=1}^{\infty}\beta^{t-1} \stageC(x^t,x^{t+1},\rxi^t)\Biggr] \leq \ctgF_{\MPC}(x^1,\mu) \quad \text{where}~x^{t+1} = y_{\MPC}(x^t,\rxi^t) \;  \forall t\in \nats.
\end{equation*}
That is, the value $\ctgF_{\MPC}(x^1,\mu)$ obtained by solving the MPC functional equation (\ref{equation:mpc-value-function}) is an upper bound on the cost incurred when applying the resulting MPC policy to (\ref{problem:soc-sum}).
\end{repeattheorem}

\proof{Proof}
Let $\mathcal{P} = \{\Prob, \pmProb_{\mu}\}$, and recall that $\mathcal{P}^\infty = \mathcal{P} \times \mathcal{P} \times\cdots$ is a rectangular set. By Theorem~\ref{theorem:mpc-is-dro}, the MPC policy solves (\ref{problem:dro-sum}), and we have that
\begin{align}
\ctgF_{\MPC}(x^1,\mu) &= \inf_{y^1,y^2, \ldots}\sup_{\Qrob\in\mathcal{P}^{\infty}} \Expt_{\Qrob}\Biggl[\,\sum_{t=1}^{\infty}\beta^{t-1}\stageC(x^t,x^{t+1},\rxi^t) \Biggr] \quad \subjectTo~x^{t+1} = y^t(x^t,\rxi^t) \;  \forall t\in \nats \notag \\
&= \sup_{\Qrob\in\mathcal{P}^{\infty}} \Expt_{\Qrob}\Biggl[\,\sum_{t=1}^{\infty}\beta^{t-1}\stageC(x^t,x^{t+1},\rxi^t) \Biggr] \quad \text{where}~x^{t+1} = y_{\MPC}(x^t,\rxi^t) \;  \forall t\in \nats. \label{equation:infsup-MPCsup}
 \end{align}
The cost incurred when applying the MPC policy to (\ref{problem:soc-sum}) is
\begin{equation*}
 \Expt_{\Prob^{\infty}}\Biggl[\, \sum_{t=1}
^{\infty}\beta^{t-1}\stageC(x^t,x^{t+1},\rxi^t) \Biggr]
\leq  \sup_{\Qrob\in\mathcal{P}^{\infty}} \Expt_{\Qrob}\Biggl[\, \sum_{t=1}^{\infty}\beta^{t-1}\stageC(x^t,x^{t+1},\rxi^t) \Biggr] \quad \text{where}~x^{t+1} = y_{\MPC}(x^t,\rxi^t) \;  \forall t\in \nats, 
\end{equation*}
which holds since $\Prob^{\infty} \in \mathcal{P}^{\infty}$. The right-hand side equals  $\ctgF_{\MPC}(x^1,\mu)$, by virtue of (\ref{equation:infsup-MPCsup}).
\Halmos
\endproof

\begin{repeatcorollary}
Suppose that the conditions in the statement of Theorem~\ref{theorem:mpc-dro-guarantee} hold, but that the MPC Bellman operator is convexity preserving rather than concavity preserving. Then
\begin{equation*}
    \Expt_{\Prob^\infty}\Biggl[\,\sum_{t=1}^{\infty}\beta^{t-1} \stageC(x^t,x^{t+1},\rxi^t)\Biggr] \geq \ctgF_{\MPC}(x^1,\mu) \quad \text{where}~x^{t+1} = y_{\MPC}(x^t,\rxi^t) \;  \forall t\in \nats;
\end{equation*}
that is, the value $\ctgF_{\MPC}(x^1,\mu)$ obtained by solving the MPC functional equation (\ref{equation:mpc-value-function}) is a lower bound on the cost incurred when applying the resulting MPC policy to (\ref{problem:soc-sum}).
\end{repeatcorollary}

\proof{Proof of Corollary~\ref{corollary:mpc-doo-dissapointment}}
The proof repeats the proof of Theorem~\ref{theorem:mpc-dro-guarantee}, but with the reference to Theorem~\ref{theorem:mpc-is-dro} replaced by a reference to Corollary~\ref{corollary:mpc-is-doo}, the supremum operations replaced by infimum operations, and the inequalities reversed. \Halmos
\endproof
\medskip

\section{Proof of Theorem~\ref{theorem:policy-switch-improvement}}
We prove the theorem by induction. To do this we extend Definition~\ref{definition:switch-policies-value-function} as follows.
\begin{definition}
For two admissible policies $y,y^{\prime}:\mathcal{X}\times\Xi\mapsto\mathcal{X}$, let $\oExptCtgF^{1}_{y,y^{\prime}}(x) \defeq \oExptCtgF_{y,y^{\prime}}(x)$, and for each $T \in \nats$ define $\oExptCtgF^{T+1}_{y,y^{\prime}}:\mathcal{X}\to\reals$ by
\begin{equation*}
\oExptCtgF^{T+1}_{y,y^{\prime}}(x) \defeq \Expt_{\Prob}\Bigl[\stageC\bigl(x,{y}(x,\rxi),\rxi\bigr)
+\beta\oExptCtgF^{T}_{y,y^{\prime}}\bigl({y}(x,\rxi)\bigr)\Bigr].
\label{equation:policy-switch-recursion}
\end{equation*}
\end{definition}
The value $\oExptCtgF^{T}_{y,y^{\prime}}(x)$ is the objective value of (\ref{problem:soc-sum}) when starting from initial state $x$ if policy ${y}$ is followed in the first $T$ stages and policy ${y^{\prime}}$ is followed thereafter. It is clear that $\oExptCtgF^{T}_{y,y^{\prime}}$ is well defined and bounded in the same way that $\oExptCtgF_{y}$ and $\oExptCtgF_{y^{\prime}}$ are.

\begin{lemma}
Let Assumption~\ref{assumption:unbounded-noise} hold. For two admissible policies ${y},{y^{\prime}}:\mathcal{X}\times\Xi\to\mathcal{X}$, it holds that $\lim_{T\to\infty}\bigl\lvert\oExptCtgF^{T}_{{y},{y^{\prime}}}(x)-\oExptCtgF_{y}(x)\bigr\rvert=0$ for every $x\in\mathcal{X}$.
\label{proposition:policy-switch-limit-is-sdp}
\end{lemma}

\proof{Proof of Lemma~\ref{proposition:policy-switch-limit-is-sdp}.}
For $x^1\in\mathcal{X}$ observe that
\begin{align*}
    \oExptCtgF^1_{{y},{y^{\prime}}}(x^1) - \oExptCtgF_{y}(x^1) &= \Expt_{\Prob}\Bigl[ \stageC\bigl(x^1,y(x^1,\rxi),\rxi\bigr)+\beta\oExptCtgF_{{y^{\prime}}}\bigl(y(x^1,\rxi)\bigr)\Bigr] - \Expt_{\Prob}\Bigl[ \stageC\bigl(x^1,y(x^1,\rxi),\rxi\bigr)+\beta\oExptCtgF_{y}\bigl(y(x^1,\rxi)\bigr)\Bigr]\\
    &= \beta\Expt_{\Prob}\Bigl[\oExptCtgF_{{y^{\prime}}}\bigl(y(x^1,\rxi)\bigr) - \oExptCtgF_{y}\bigl(y(x^1,\rxi)\bigr)\Bigr].
\end{align*}
Continuing this reasoning, an induction shows that 
\begin{align}\label{equation:difference-of-switch-policy-oos-value}
    \bigl\lvert \oExptCtgF^T_{{y},{y^{\prime}}}(x^1) - \oExptCtgF_{y}(x^1)\bigr\rvert = \beta^T \Bigl\lvert \Expt_{\Prob^\infty}\bigl[\oExptCtgF_{{y^{\prime}}}(x^T) - \oExptCtgF_{y}(x^T)\bigr]\Bigr\rvert \quad \text{where}~x^{t+1} = y(x^t,\rxi^t) \; \forall t\in\{1,\ldots,T-1\}.
\end{align}
Since $\oExptCtgF_{y}$ and $\oExptCtgF_{{y^{\prime}}}$ are bounded under Assumption~\ref{assumption:unbounded-noise} and $\beta^T\to 0$ as $T\to\infty$, the value of (\ref{equation:difference-of-switch-policy-oos-value}) vanishes as $T\rightarrow \infty$. Replacing $x^1$ with $x$ concludes the proof.
\Halmos
\endproof

\begin{repeattheorem}[Theorem~\ref{theorem:policy-switch-improvement}.]
Let Assumption~\ref{assumption:unbounded-noise} hold. For two admissible policies $y,y^{\prime}:\mathcal{X}\times\Xi\to\mathcal{X}$, if
\begin{equation*}
     \oExptCtgF_{y^{\prime}}(x)\leq\oExptCtgF_{y,y^{\prime}}(x) \quad \forall x \in \mathcal{X},
\end{equation*}
then $\oExptCtgF_{y^{\prime}}(x)\leq\oExptCtgF_{y}(x)$ for each $x \in \mathcal{X}$.
\end{repeattheorem}

\proof{Proof of Theorem~\ref{theorem:policy-switch-improvement}.}
Denoting $\oExptCtgF_{{y^{\prime}}}(x)$ by $\oExptCtgF^{0}_{{y},{y^{\prime}}}(x)$, we pose the inductive hypothesis $\oExptCtgF^{t-1}_{{y},{y^{\prime}}}(x)\leq \oExptCtgF^{t}_{{y},{y^{\prime}}}(x)$ for all $x\in\mathcal{X}$, which is true for $t=1$ by assumption. By the hypothesis, ${\oExptCtgF^{t-1}_{{y},{y^{\prime}}}({y}(x,\xi)) \leq \oExptCtgF^{t}_{{y},{y^{\prime}}}(y(x,\xi))}$, and it follows that
\begin{align}
\oExptCtgF^{t}_{{y},{y^{\prime}}}(x) & =\Expt_{\Prob}\Bigl[\stageC\bigl(x,{y}(x,\rxi),\rxi\bigr)
+\beta \oExptCtgF^{t-1}_{{y},{y^{\prime}}}\bigl({y}(x,\rxi)\bigr)\Bigr]   \notag \\
& \leq \Expt_{\Prob}\Bigl[ \stageC\bigl(x,{y}(x,\rxi),\rxi\bigr)
+\beta\oExptCtgF^{t}_{{y},{y^{\prime}}}\bigl({y}(x,\rxi)\bigr)\Bigr] 
 = \oExptCtgF^{t+1}_{{y},{y^{\prime}}}(x)
\label{equation:policy-switch-induction-base-case}
\end{align}
for all $x\in\mathcal{X}$. The inequality (\ref{equation:policy-switch-induction-base-case}) establishes the induction. Hence, $\oExptCtgF_{{y^{\prime}}}(x) \leq \oExptCtgF^{t}_{{y},{y^{\prime}}}(x)$ for each $t\in\nats$ and ${\oExptCtgF_{{y^{\prime}}}(x) \leq \lim_{t\rightarrow \infty }\oExptCtgF^{t}_{{y},{y^{\prime}}}(x) = \oExptCtgF_{y}(x)}$, where Lemma~\ref{proposition:policy-switch-limit-is-sdp} yields the final equality.\Halmos
\endproof
\medskip

\section{Proof of Proposition~\ref{proposition:rosp-optimal-policy}}
\begin{repeatproposition}[Proposition~\ref{proposition:rosp-optimal-policy}.]
Suppose that the distribution $\Prob$ has compact support. Then the policy which at inventory $x$ and price $p$ sells down to
\begin{equation*}
y(x,p) = \storageCostDrv^{-1}\Bigl(\bigl( \beta  \Expt_{\Prob}\bigl[(\rp-{p})_+\bigr]-(1-\beta)p\bigr)_{[\storageCostDrv(0),\storageCostDrv(x)]} \Bigr)
\end{equation*}
is optimal for (\ref{problem:rosp-sum}).
\end{repeatproposition}

\proof{Proof of Proposition~\ref{proposition:rosp-optimal-policy}.}
Recall that there is a unique continuous function $\ctgF:\reals_+\times\reals\to\reals$ satisfying  (\ref{equation:rosp-recursion}) and that $x\mapsto v(x,p)$ is concave for each $p\in\reals$. Given such $v$ we derive a policy which satisfies (\ref{equation:rosp-argmax}) and therefore solves (\ref{problem:rosp-sum}). Define $\ExptCtgF:\reals_+\to\reals$ by $\ExptCtgF(x) \defeq \Expt_{\Prob}\bigl[\ctgF(x,\rp)\bigr]$ and note that $V$ is continuous and concave \citep[Theorems~7.43, 7.46]{shapiro-et-al:lectures-on-SP}. Satisfying (\ref{equation:rosp-argmax}) amounts to solving 
\begin{equation}\label{problem:optimal-policy}
\maximize_{y\in\reals_+} ~~ p(x-y)-\storageCost(y)+\beta \ExptCtgF(y)~~\subjectTo~~ y \in[0,x].
\end{equation}
Due to the strict concavity of $C$, (\ref{problem:optimal-policy}) has a unique solution $y^{\star}(x,p)\in [0,x]$. Let $y^{\star}(\infty,p) \defeq \lim_{x\to\infty} y^{\star}(x,p)$. (This limit is finite valued since $\lim_{y\to\infty}-c(y)=-\infty$.) Now, $y^{\star}(\infty,p)$ is equal to the optimal solution $y^{\star}(x,p)$ when projected onto $[0,x]$, and the mapping $p\mapsto y^{\star}(\infty,p)$ is decreasing. Thus there exists $\cPrice(x)\in[-\infty,\infty]$ such that for each $p\in\reals$, if $p\leq\cPrice(x)$, then $y^{\star}(\infty,p)\geq x$, and otherwise if $p\geq\cPrice(x)$, then $y^{\star}(\infty,p) \leq x$. Comparing $\cPrice(x)$ and $p$ characterises whether or not $y^{\star}(\infty,p)$ needs to be projected onto $[0,x]$ to solve (\ref{problem:optimal-policy}).

We proceed by deriving an expression for the derivative of $\ExptCtgF$ using a fixed-point equation. In fact we use the superdifferential operator $\partial$ since $\ExptCtgF$ may be nonsmooth. Due to continuity and concavity, the superdifferential $\partial\ExptCtgF(x)$ is nonempty and compact. We also denote the superdifferential of $\ctgF$ with respect to $x$ by $\partial_x \ctgF(x,p)$. When $p \leq \cPrice(x)$, by definition $y^{\star}(\infty,p) \geq x$, and thus $y^{\star}(x,p)=x$. In view of (\ref{equation:rosp-recursion}),
\begin{align}\label{equation:-cost+beta-times-value-subgrad}
\ctgF(x,p)=-\storageCost(x)+\beta \ExptCtgF(x)& \quad \forall p\leq \cPrice(x)\notag\\
\implies -\storageCostDrv(x)+\beta \partial{\ExptCtgF}(x)\subseteq \partial_x \ctgF(x,p)&.
\end{align}
On the other hand, 
\begin{align}\label{equation:p-subgrad}
    \ctgF(x,p)=p\bigl(x-y^{\star}(\infty,p)\bigr)-\storageCost\bigl(y^{\star}(\infty,p)\bigr)+\beta \ExptCtgF\bigl(y^{\star}(\infty,p)\bigr)& \quad \forall p\geq \cPrice(x)\notag\\
    \implies p \in \partial_x \ctgF(x,p)&.
\end{align}
We therefore define $\psi_x:\partial{\ExptCtgF}(x)\times\reals\to\reals$ by
\begin{equation}\label{equation:piecewise-gradient}
\psi_x(g,p) \defeq
\begin{cases}
-\storageCostDrv(x)+\beta g & \text{if}~p \leq \cPrice(x)  \\ 
p & \text{if}~p > \cPrice(x);
\end{cases}
\end{equation}
the two cases here coinciding with (\ref{equation:-cost+beta-times-value-subgrad}) and (\ref{equation:p-subgrad}). Hence, for each $g\in\partial{\ExptCtgF}(x)$,
\begin{equation*}
    \psi_x(g,p) \in \partial_x \ctgF(x,p) \quad \forall p\in\reals,
\end{equation*}
and thus $\Expt_{\Prob}\bigl[\psi_x(g,\rp)\bigr] \in \partial {\ExptCtgF}(x)$ \citep[Theorem~7.47]{shapiro-et-al:lectures-on-SP}. The expectation of (\ref{equation:piecewise-gradient}) is
\begin{equation}\label{equation:expected-subgrad}
\Expt_{\Prob}\bigl[\psi_x(g,\rp)\bigr] = \Expt_{\Prob}\bigl[\Indicator\{\rp\leq\cPrice(x)\}\bigr]\bigl(\beta g-\storageCostDrv(x)\bigr)+\Expt_{\Prob}\bigl[\rp\Indicator\{\rp > \cPrice(x)\}\bigr],
\end{equation}
and for $g,g^{\prime}\in\partial{V}(x)$ taking the difference shows that
\begin{equation*}
\Bigl\lvert \Expt_{\Prob}\bigl[\psi_x(g,\rp)\bigr] -\Expt_{\Prob}\bigl[\psi_x(g^{\prime},\rp)\bigr]\Bigr\rvert = \Expt_{\Prob}\bigl[\Indicator\{\rp \leq \cPrice(x)\}\bigr]\beta\lvert g-g^{\prime}\rvert < \lvert g-g^{\prime}\rvert,
\end{equation*}
i.e., $g\mapsto\Expt_{\Prob}\bigl[\psi_x(g,\rp)\bigr]$ is a contraction mapping. With $(\partial\ExptCtgF(x),\lvert\blank\rvert)$ a complete metric space, by Banach's fixed-point theorem, there is a unique ${g(x) \in \partial\ExptCtgF(x)}$ satisfying $g(x) = \Expt_{\Prob}\bigl[\psi_x(g(x),\rp)\bigr]$. So (\ref{equation:expected-subgrad}) yields
\begin{align}
g(x) &= \Expt_{\Prob}\bigl[\Indicator\{\rp\leq \cPrice(x)\}\bigr]\bigl(\beta g(x)-\storageCostDrv(x)\bigr)+\Expt_{\Prob}\bigl[\rp\Indicator\{\rp > \cPrice(x)\}\bigr]\notag\\
&= \frac{\Expt_{\Prob}\bigl[\rp\Indicator\{\rp > \cPrice(x)\}\bigr]-\Expt_{\Prob}\bigl[\Indicator\{\rp \leq \cPrice(x)\}\bigr]\storageCostDrv(x)}{1-\beta\Expt_{\Prob}\bigl[\Indicator\{\rp \leq \cPrice(x)\}\bigr]}\in\partial\ExptCtgF(x).
\end{align}
Using this supergradient and observing that $p$ is by definition the minimum acceptable price $p^{\star}(y(\infty,p))$, a first-order optimality condition for the unconstrained version of (\ref{problem:optimal-policy}) is
\begin{align*}
0 &= -p-c(y) +\beta\frac{\Expt_{\Prob}\bigl[\rp\Indicator\{\rp > p\}\bigr]-\Expt_{\Prob}\bigl[\Indicator\{\rp \leq p\}\bigr]\storageCostDrv(y)}{1-\beta\Expt_{\Prob}\bigl[\{\rp \leq p\}\bigr]}\notag\\
&= -\bigl( 1-\beta \Expt_{\Prob}\bigl[ \Indicator\{\rp \leq p\}\bigr] \bigr)p -c(y) +\beta\Expt_{\Prob}\bigl[\rp\Indicator\{\rp > p\}\bigr],
\end{align*}
which can be rearranged to
\begin{align}\label{equation:policy-first-order-condition}
c(y) &= -\bigl( 1-\beta \Expt_{\Prob}\bigl[ \Indicator\{\rp \leq p\}\bigr]\bigr)p +\beta \Expt_{\Prob}\bigl[\rp\Indicator\{\rp > p\}\bigr] \notag\\
&= \beta \Expt_{\Prob}\bigl[\rp\Indicator\{\rp > p\}\bigr] - \beta \Expt_{\Prob}\bigl[ p\Indicator\{\rp > p\}\bigr]  -\bigl( 1-\beta \bigr)p  \notag\\
&= \beta \Expt_{\Prob}\bigl[( \rp-p) _{+}\bigr] -( 1-\beta) p.
\end{align}
If this is an element of $[\storageCostDrv(0),\storageCostDrv(x)]$, then $\storageCostDrv^{-1}\bigl( \beta\Expt_{\Prob}\bigl[(\rp-p)_{+}\bigr]-( 1-\beta) p\bigr)$ solves (\ref{problem:optimal-policy}). Otherwise, due to the concavity of the univariate objective function, first projecting (\ref{equation:policy-first-order-condition}) onto $[\storageCostDrv(0),\storageCostDrv(x)]$ solves (\ref{problem:optimal-policy}).
\Halmos
\endproof

\section{Proof of Proposition~\ref{proposition:MPC-wins-with-outliers}}
We prove the result using several lemmas which provide formulae for $\frac{\drv}{\drv x}\oExptCtgF_{\SDP}(x)$ and $\frac{\drv}{\drv x}\oExptCtgF_{\MPC}(x)$. To this end the following lemma will be useful.
\begin{lemma}{\citep[Theorem~7.44]{shapiro-et-al:lectures-on-SP}}\label{lemma:Shapiro-derivatives}
Let $f:\reals^m \times \Omega \rightarrow \reals$ be a random function with expected value $F(x)\defeq\Expt\bigl[f(x,\romega)\bigr]$. If the following conditions hold:
\begin{enumerate}[label={\normalfont(\roman*)}, ref={\roman*}]
\item\label{lemma:Shapiro-derivatives-i} The expectation $F(\bar{x})$ is well defined and finite valued at a given point $\bar{x} \in \reals^m$,
\item\label{lemma:Shapiro-derivatives-ii} There exists a positive-valued random variable $l(\romega)$ such that $\Expt\bigl[l(\romega)\bigr] < \infty$, and for each $x,x^{\prime}$ in a neighbourhood of $\bar{x}$ and almost every $\omega\in\Omega$ the following inequality holds:
\begin{equation*}
    \bigl\lvert f(x,\omega)-f(x^{\prime},\omega)\bigr\rvert \leq l(\omega) \lVert x - x^{\prime} \rVert_2,
\end{equation*}
\item\label{lemma:Shapiro-derivatives-iii} For almost every $\omega\in\Omega$ the function $f(x,\omega)$ is differentiable with respect to $x$ at $\bar{x}$,
\end{enumerate}
Then $F(x)$ is differentiable at $\bar{x}$ and
\begin{equation*}
     \nabla F(\bar{x})= \Expt\bigl[\nabla_x f(\bar{x},\romega)\bigr].
 \end{equation*}
\end{lemma}

Now the derivative values can be established.
\begin{lemma}
Let Assumption~\ref{assumption:rosp-unbounded-prices} hold. Then the functions $\oExptCtgF_{\SDP}$ and $\oExptCtgF_{\MPC}$ are differentiable.
\label{lemma:oos-derivatives-exist}
\end{lemma}

\proof{Proof of Lemma~\ref{lemma:oos-derivatives-exist}.}
We proceed by verifying that the conditions (\ref{lemma:Shapiro-derivatives-i})--(\ref{lemma:Shapiro-derivatives-iii}) of Lemma~\ref{lemma:Shapiro-derivatives} apply to $\oExptCtgF_{\SDP}$. Under Assumption~\ref{assumption:rosp-unbounded-prices}, note that $\oExptCtgF_{\SDP}$ is well defined and finite valued, satisfying (\ref{lemma:Shapiro-derivatives-i}). Next we address (\ref{lemma:Shapiro-derivatives-iii}). Define $\oCtgF_{\SDP}:\reals_+\times\reals_+^{\infty}\to\reals$ by
\begin{equation}\label{equation:instance-of-sum-of-stage-reward}
\oCtgF_{\SDP}(x^{1},p^{1},p^{2},\ldots) \defeq \sum_{t=1}^{\infty}\beta^{t-1}\bigl(p^{t}(x^{t}-x^{t+1})-\storageCost(x^{t+1})\bigr) \quad \text{where}~x^{t+1} = y_{\SDP}(x^{t},p^{t}) \; \forall t\in\nats. 
\end{equation}
If $\rp^1,\rp^2,\ldots$ are random price variables that are independent and identically distributed according to $\Prob$, the definition (\ref{equation:out-of-sample-sum}) shows that $\Expt_{\Prob^\infty}\bigl[\oCtgF_{\SDP}(x^1,\rp^1,\rp^2,\ldots)\bigr] = \oExptCtgF_{\SDP}(x^1)$. Recall that $\cPriceSDP(x^1)$ is the minimum price required for the SDP policy to decide that it is worthwhile to sell a portion of inventory from $x^1$, and that for any $p \geq \cPriceSDP(x^1)$ the SDP policy sells down to the target inventory level $y_{\SDP}(\infty,p) \leq x^1$. It follows that in (\ref{equation:instance-of-sum-of-stage-reward}) all inventory levels $x^t = x^1$ for $t \leq T$, with $T$ being the first time at which $p^T \geq \cPriceSDP(x^1)$. The next inventory level $x^{T+1} = y_{\SDP}(\infty,p^T)$, and along with the prices $p^{T+1},p^{T+2},\ldots,$ these uniquely determine the remaining inventory levels $x^{T+2},x^{T+3},\ldots$ when following the SDP policy. Under Assumption~\ref{assumption:rosp-unbounded-prices} the probability distribution $\Prob$ is atomless which implies that $\Prob\bigl[\rp^t \neq \cPriceSDP(x^{1})\bigr]=1$. It follows that in (\ref{equation:instance-of-sum-of-stage-reward}) there are no $p^{t}$ with $y_{\SDP}(\infty,p^{t})=x^{1}$ almost surely. Hence, there is a neighbourhood of $x^{1}$ values in which $x^t=x^1$ for $t \leq T$ and $x^{T+1} = y_{\SDP}(\infty,p^T)$. The remaining inventory values $x^{T+2},x^{T+3},\ldots$ are fixed within this neighbourhood and we write
\begin{equation}\label{equation:instance-derivative-of-sum-of-stage-reward}
    \frac{\drv \oCtgF_{\SDP}(x^{1},p^{1},p^{2},\ldots)}{\drv x^1}= \beta^{T-1}p^{T} - \sum_{t=1}^{T-1}\beta^{t-1} \storageCostDrv(x^1),
\end{equation}
meeting (\ref{lemma:Shapiro-derivatives-iii}). 

It remains to verify (\ref{lemma:Shapiro-derivatives-ii}). When attempting to evaluate $\Expt_{\Prob^\infty}\bigl[\frac{\drv}{\drv x^1}\oCtgF_{\SDP}(x^{1},\rp^{1},\rp^{2},\ldots)\bigr]$ using the expression (\ref{equation:instance-derivative-of-sum-of-stage-reward}), the time $\rT$ and the price $\rp^{\rT}$ are random variables. The expectation of the first term in (\ref{equation:instance-derivative-of-sum-of-stage-reward}) is 
\begin{equation*}
\Expt_{\Prob^\infty}\bigl[\beta^{\rT-1} \rp^{\rT}\bigr]\leq\Expt_{\Prob}\bigl[\rp\givenn\rp \geq \cPriceSDP(x^1)\bigr] \leq \Expt_{\Prob}\bigl[\rp\givenn\rp\geq \max\{{p}_1,\ldots,{p}_N\}\bigr] < \infty,
\end{equation*}
where the first inequality follows since $\rp^{\rT}$ is the first price greater than $\cPriceSDP(x^1)$, the second inequality follows since the SDP policy always decides that it is worthwhile to sell when the price is higher than the highest price sample, and the third inequality follows since $\Prob$ has a finite mean. The expectation of the second term in (\ref{equation:instance-derivative-of-sum-of-stage-reward}) is $\Expt_{\Prob^\infty}\bigl[\sum_{t=1}^{\rT-1}\beta^{t-1} \storageCostDrv(x^1)\bigr] \leq \frac{1}{1-\beta}\storageCostDrv(x^1)$ which is bounded on compact sets. Together these observations show that $\oCtgF_{\SDP}(x^{1},p^{1},p^{2},\ldots)$ has a Lipschitz constant in its first argument with finite expectation, meeting (\ref{lemma:Shapiro-derivatives-ii}). Thus, Lemma~\ref{lemma:Shapiro-derivatives} applies and $\oExptCtgF_{\SDP}(x)$ is differentiable. Similar reasoning shows $\oExptCtgF_{\MPC}(x)$ is differentiable. \Halmos
\endproof

\begin{lemma}
Let Assumption~\ref{assumption:rosp-unbounded-prices} hold. Then the derivatives of the functions $\oExptCtgF_{\SDP}(x)$ and $\oExptCtgF_{\MPC}(x)$ are given by
\begin{equation*}
\frac{\drv\oExptCtgF_{\SDP}(x)}{\drv x} = \frac{\Expt_{\Prob}\bigl[\rp\Indicator\{\rp> \cPriceSDP(x)\}\bigr]-\Expt_{\Prob}\bigl[\Indicator\{\rp\leq\cPriceSDP(x)\}\bigr]\storageCostDrv(x)}{1-\beta \Expt_{\Prob}\bigl[\Indicator\{\rp \leq \cPriceSDP(x)\}\bigr]}
\end{equation*}
and 
\begin{equation*}
\frac{\drv\oExptCtgF_{\MPC}(x)}{\drv x}=
\frac{\Expt_{\Prob}\bigl[\rp\Indicator\{\rp> \cPriceMPC(x)\}\bigr]-\Expt_{\Prob}\bigl[\Indicator\{\rp\leq\cPriceMPC(x)\}\bigr]\storageCostDrv(x)}{1-\beta \Expt_{\Prob}\bigl[\Indicator\{\rp \leq \cPriceMPC(x)\}\bigr]},
\end{equation*}
respectively.
\label{lemma:oos-derivatives}
\end{lemma}

\proof{Proof of Lemma~\ref{lemma:oos-derivatives}.}
We proceed by showing that the derivative $\frac{\drv}{\drv x}\oExptCtgF_{\SDP}(x)$ satisfies a functional equation. Define $\oCtgF_{\SDP}:\reals_+\times\reals_+\to\reals$ by
\begin{equation}
\oCtgF_{\SDP}(x,p) \defeq
        \begin{cases} 
            -\storageCost(x)+\beta \oExptCtgF_{\SDP}(x) &  \text{if}~ p\leq {\cPriceSDP(x)}\\
            p\bigl(x-y_{\SDP}(\infty,p)\bigr)-\storageCost\bigl(y_{\SDP}(\infty,p)\bigr)+\beta \oExptCtgF_{\SDP}\bigl(y_{\SDP}(\infty,p)\bigr) & \text{if}~ {p} > {\cPriceSDP(x)}.
        \end{cases}
\label{equation:piece-wise-out-of-sample-value}
\end{equation}
The equation (\ref{equation:rosp-oos-sdp-value-function}) shows $\Expt_{\Prob}\bigl[ \oCtgF_{\SDP}(x,\rp)\bigr]=\oExptCtgF_{\SDP}(x)$. Recalling that $\cPriceSDP$ is continuous, for each $p<\cPriceSDP(x)$ there is a neighbourhood of $x$ values for which the first case of (\ref{equation:piece-wise-out-of-sample-value}) holds. Due to Lemma~\ref{lemma:oos-derivatives-exist}, $\frac{\drv}{\drv x}\oExptCtgF_{\SDP}(x)$ exists, and thus
\begin{equation}\label{equation:piece-wise-out-of-sample-derivative-A}
    \frac{\drv \oCtgF_{\SDP}(x,p)}{\drv x} = -\storageCostDrv(x)+\beta\frac{\drv \oExptCtgF_{\SDP}(x)}{\drv x}.
\end{equation}
In the other case, for each $p>\cPriceSDP(x)$ the second case of (\ref{equation:piece-wise-out-of-sample-value}) holds and 
\begin{equation}\label{equation:piece-wise-out-of-sample-derivative-B}
    \frac{\drv \oCtgF_{\SDP}(x,p)}{\drv x} = p.
\end{equation}
Similar reasoning to that in the proof of Lemma~$\ref{lemma:oos-derivatives-exist}$ shows $\frac{\drv}{\drv x}\oExptCtgF_{\SDP}(x)=\Expt_{\Prob}\bigl[\frac{\drv}{\drv x}\oCtgF_{\SDP}(x,\rp)\bigr]$, and with $\Prob$ atomless, the event $\rp = \cPriceSDP(x)$ is immaterial when evaluating $\Expt_{\Prob}\bigl[\frac{\drv}{\drv x}\oCtgF_{\SDP}(x,\rp)\bigr]$. Taking the expectation of (\ref{equation:piece-wise-out-of-sample-derivative-A}) and (\ref{equation:piece-wise-out-of-sample-derivative-B}) yields
\begin{align*}
\frac{\drv \oExptCtgF_{\SDP}(x)}{\drv x} &=\Expt_{\Prob}\Bigl[\Indicator\bigl\{\rp \leq \cPriceSDP(x)\bigr\}\Bigr]\biggl( \beta \frac{\drv \oExptCtgF_{\SDP}(x)}{\drv x} -\storageCostDrv(x)\biggr) + \Expt_{\Prob}\Bigl[\rp\Indicator\bigl\{\rp> \cPriceSDP(x)\bigr\}\Bigr]\\
&=\frac{\Expt_{\Prob}\bigl[\rp\Indicator\{\rp> \cPriceSDP(x)\}\bigr]-\Expt_{\Prob}\bigl[\Indicator\{\rp\leq\cPriceSDP(x)\}\bigr]\storageCostDrv(x)}{1-\beta \Expt_{\Prob}\bigl[\Indicator\{\rp \leq \cPriceSDP(x)\}\bigr]}.
\end{align*}
The expression for $\frac{\drv}{\drv x}\oExptCtgF_{\MPC}(x)$ can be derived in the same way.
\Halmos \endproof

Now Proposition~\ref{proposition:MPC-wins-with-outliers} can be established.
\begin{repeatproposition}[Proposition~\ref{proposition:MPC-wins-with-outliers}.]
Let Assumption~\ref{assumption:rosp-unbounded-prices} hold. For price samples $p_1,\ldots,p_N$ which determine the sample-based SDP and MPC policies by (\ref{equation:sdp-policy}) and (\ref{equation:mpc-policy}), if the minimum acceptable SDP price $\cPriceSDP(x;p_1,\ldots,p_N)$ is such that
\begin{equation*}
    \storageCostDrv(x) \geq \beta \Expt_{\Prob}\Bigl[\,\rp\Indicator\bigl\{\rp > \cPriceSDP(x;p_1,\ldots,p_N)\bigr\}\Bigr]  \quad \forall x\in [0,x^1], 
\end{equation*}
then $\oExptCtgF_{\MPC}(x;p_1,\ldots,p_N) \geq \oExptCtgF_{\SDP}(x;p_1,\ldots,p_N)$ for each $x\in [0,x^1]$. That is, MPC performs at least as well out of sample on (\ref{problem:rosp-sum}) as SDP.
\end{repeatproposition}

\proof{Proof of Proposition~\ref{proposition:MPC-wins-with-outliers}.}
Starting from initial inventory $x\in[0,x^1]$, let $\oExptCtgF_{\SDP,\MPC}(x)$ denote the objective value of (\ref{problem:rosp-sum}) when following the SDP policy in the first stage and the MPC policy thereafter (as for Definition~\ref{definition:switch-policies-value-function}). We first compare derivatives and show ${\frac{\drv}{\drv x}\oExptCtgF_{\MPC}(x)\geq \frac{\drv}{\drv x}\oExptCtgF_{\SDP,\MPC}(x)}$.
Rearranging the identity provided by Lemma~\ref{lemma:oos-derivatives},
\begin{equation}\label{equation:MPC-value-derivative}
\frac{\drv\oExptCtgF_{\MPC}(x)}{\drv x}= \Expt_{\Prob}\Bigl[\Indicator\bigl\{\rp\leq\cPriceMPC(x)\bigr\}\Bigr]\biggl(\beta \frac{\drv\oExptCtgF_{\MPC}(x)}{\drv x}-\storageCostDrv(x)\biggr) + \Expt_{\Prob}\Bigl[\rp\Indicator\bigl\{\rp>\cPriceMPC(x)\bigr\}\Bigr],
\end{equation}
and similarly
\begin{equation*}
\frac{\drv \oExptCtgF_{\SDP,\MPC}(x)}{\drv x} = \Expt_{\Prob}\Bigl[\Indicator\bigl\{\rp\leq\cPriceSDP(x)\bigr\}\Bigr]\biggl(\beta \frac{\drv\oExptCtgF_{\MPC}(x)}{\drv x}-\storageCostDrv(x)\biggr) + \Expt_{\Prob}\Bigl[\rp\Indicator\bigl\{\rp>\cPriceSDP(x)\bigr\}\Bigr].
\end{equation*}
Thus, to compare the derivatives it suffices to show
\begin{equation}\label{equation:MPC-wins-with-outliers}
\Expt_{\Prob}\Bigl[\rp\Indicator\bigl\{\cPriceMPC(x) < \rp\leq\cPriceSDP(x)\bigr\}\Bigr] \geq \Expt_{\Prob}\Bigl[\Indicator\bigl\{\cPriceMPC(x) < \rp\leq\cPriceSDP(x)\bigr\}\Bigr]\biggl( \beta \frac{\drv \oExptCtgF_\MPC(x)}{\drv x} - \storageCostDrv(x) \biggr).
\end{equation}
From (\ref{equation:MPC-value-derivative}) we have
\begin{equation*}\label{equation:MPC-value-derivative-inequality}
\beta \frac{\drv \oExptCtgF_\MPC(x)}{\drv x} - \storageCostDrv(x) = \frac{ \beta \Expt_{\Prob}\bigl[\rp\Indicator\{\rp > \cPriceMPC(x)\}\bigr] - \storageCostDrv(x) }{1 - \beta \Expt_{\Prob}\bigl[\Indicator\{\rp \leq \cPriceMPC(x)\}\bigr]}
\leq \frac{ \beta \Expt_{\Prob}\bigl[\rp\Indicator\{ \cPriceMPC(x) < \rp \leq \cPriceSDP(x)\}\bigr] }{1 - \beta \Expt_{\Prob}\bigl[\Indicator\{\rp \leq \cPriceMPC(x)\}\bigr]},
\end{equation*}
where the inequality follows from the condition in the statement of the proposition. With ${\beta \in (0,1)}$ and $1 - \beta \Expt_{\Prob}\bigl[\Indicator\{\rp \leq \cPriceMPC(x)\}\bigr] \geq \Expt_{\Prob}\bigl[\Indicator\{\cPriceMPC(x) < \rp \leq \cPriceSDP(x)\}\bigr]$, multiplying by $\Expt_{\Prob}\bigl[\Indicator\{\cPriceMPC(x) < \rp\leq\cPriceSDP(x)\}\bigr]$ yields (\ref{equation:MPC-wins-with-outliers}), whereby ${\frac{\drv}{\drv x}\oExptCtgF_{\MPC}(x)}\geq {\frac{\drv}{\drv x}\oExptCtgF_{\SDP,\MPC}(x)}$. Hence, ${\oExptCtgF_{\MPC}(x)\geq \oExptCtgF_{\SDP,\MPC}(x)}$ for all $x\in [0,x^1]$, and from Theorem~\ref{theorem:policy-switch-improvement} we conclude $\oExptCtgF_{\MPC}(x)\geq \oExptCtgF_{\SDP}(x)$ for each $x\in [0,x^1]$, as required.
\Halmos \endproof
\medskip

\section{Proof of Proposition~\ref{proposition:exponential-example}}
To make the sample-dependence explicit, we write $\oExptCtgF_{\SDP}(x;p_1,\ldots,p_N) = \oExptCtgF_{\SDP}(x)$ and $ \oExptCtgF_{\MPC}(x;p_1,\ldots,p_N) = \oExptCtgF_{\MPC}(x)$. The expectations $\Expt_{\Prob^N}\bigl[\oExptCtgF_{\SDP}(x;\rp_1,\ldots,\rp_N)\bigr]$ and $\Expt_{\Prob^N}\bigl[\oExptCtgF_{\MPC}(x;\rp_1,\ldots,\rp_N)\bigr]$ are then the expected out-of-sample performances of the sample-based policies.

\begin{repeatproposition}[Proposition~\ref{proposition:exponential-example}]
    Let $\Prob$ be the $\text{\normalfont Exponential}(1)$ distribution, $x^1 = 1$, and $\storageCost(x) = \frac{1}{2}x^{2}$. For each sample size $N \geq 2$, as the discount factor $\beta\to 1$ the expected out-of-sample performance of SDP on (\ref{problem:rosp-sum}) is unbounded below, while that of MPC is bounded.  
\end{repeatproposition}

\proof{Proof of Proposition~\ref{proposition:exponential-example}}
First consider the term $\Expt_{\Prob^N}\bigl[\oExptCtgF_{\SDP}(1; \rp_1,\ldots,\rp_N)\bigr]$. With $\oExptCtgF_{\SDP}(0;p_1,\ldots,p_N)=0$ we write $\oExptCtgF_{\SDP}(1;p_1,\ldots,p_N)=\int_{0}^{1}\frac{\drv}{\drv x}\oExptCtgF_{\SDP}(x;p_1,\ldots,p_N)\,\drv x$ using the expression for the derivative given by Lemma~\ref{lemma:oos-derivatives}. Hence,
\begin{equation}\label{equation:exponential-sdp-expected-oos-performance}
\Expt_{\Prob^N}\bigl[\oExptCtgF_{\SDP}(1;\rp_1,\ldots,\rp_N)\bigr] = \Expt_{\Prob^N}\Biggl[\int_{0}^{1}{\frac{{\int_{\cPriceSDP(x;\rp_1,\ldots,\rp_N)}^{\infty}pf(p)\,\drv p-xF(\cPriceSDP(x;\rp_1,\ldots,\rp_N))}}{1-\beta F(\cPriceSDP(x;\rp_1,\ldots,\rp_N))}}\,\drv x\Biggr],
\end{equation}
where $F$ and $f$ denote the cumulative distribution function and probability density function of the $\text{Exponential(1)}$ distribution, respectively. The negative term in (\ref{equation:exponential-sdp-expected-oos-performance}) is
\begin{multline}\label{equation:exponential-sdp-expected-oos-performance-negative-part}
\Expt_{\Prob^N}\Biggl[\int_{0}^{1}{\frac{-xF(\cPriceSDP(x;\rp_1,\ldots,\rp_N))}{1-\beta F(\cPriceSDP(x;\rp_1,\ldots,\rp_N))}}\,\drv x\Biggr] = \\
\int_{0}^{\infty} \cdots \int_{0}^{\infty} \int_{0}^{1} \frac{-xF(\cPriceSDP(x;p_1,\ldots,p_N))}{1-\beta F(\cPriceSDP(x;p_1,\ldots,p_N))}\, \drv x \,f(p_{N})\, \drv p_{N} \cdots  f(p_{1})\, \drv p_{1}.
\end{multline}
Fix $p_{1},\ldots ,p_{N-1}$ and consider the inner-most integral in (\ref{equation:exponential-sdp-expected-oos-performance-negative-part}) when $p_{N}$ is large. Following (\ref{equation:sdp-policy}), for $\alpha\in(0,1)$ and all $\beta \in [\alpha,1)$, $x\in [0,1]$, the value $\cPriceSDP(x;p_{1},\ldots ,p_{N})$ is the $p$ which solves ${\beta\frac{1}{N}(p_{N}-p)_{+}+(1-\beta)p}={x}$. Thus
\begin{equation*}
\frac{-xF(\cPriceSDP(x;p_{1},\ldots ,p_{N}))}{1-\beta F(\cPriceSDP(x;p_{1},\ldots ,p_{N}))} \rightarrow \frac{-x( 1-e^{-(p_{N}-Nx)}) }{1-( 1-e^{-(p_{N}-Nx)}) }
\end{equation*}
uniformly in $x\in [0,1]$ as $\beta\rightarrow 1$, whereby
\begin{align*}
\lim _{\beta \rightarrow 1}\int_{0}^{1}\frac{-xF(\cPriceSDP(x;p_{1},\ldots ,p_{N}))}{1-\beta F(\cPriceSDP(x;p_{1},\ldots ,p_{N}))}\, \drv x& =\int_{0}^{1}\frac{-x( 1-e^{-(p_{N}-Nx)}) }{1-( 1-e^{-(p_{N}-Nx)}) }\, \drv x =\frac{1}{2}+\biggl( \frac{1}{N}+\frac{1}{N^{2}}\biggr) e^{p_{N}-N}-\frac{1}{N^{2}}e^{p_{N}}.
\end{align*}
So for each $p_{N}$ sufficiently large, given $\varepsilon >0$ there exists a $\beta<1$ beyond which 
\begin{equation}\label{equation:exponential-inequality-beta}
\int_{0}^{1}\frac{-xF(\cPriceSDP(x;p_{1},\ldots ,p_{N}))}{1-\beta F(\cPriceSDP(x;p_{1},\ldots ,p_{N}))}\, \drv x<\frac{1}{2}+\biggl( \frac{1}{N}+\frac{1}{N^{2}}\biggr) e^{p_{N}-N}-\frac{1}{N^{2}}e^{p_{N}}+\varepsilon.
\end{equation}
But
\begin{align*}
    \lim_{p_N\to\infty}\biggl(\frac{1}{2}+\biggl( \frac{1}{N}+\frac{1}{N^{2}}\biggr) e^{p_{N}-N}-\frac{1}{N^{2}}e^{p_{N}}+\varepsilon\biggr)e^{-p_N} = \biggl( \frac{1}{N}+\frac{1}{N^{2}}\biggr) e^{-N}-\frac{1}{N^{2}}<0
\end{align*}
and for $0<a<b$ the integral
\begin{equation*}
    \int_{a}^{b}\biggl(\frac{1}{2}+\biggl( \frac{1}{N}+\frac{1}{N^{2}}\biggr) e^{p_{N}-N}-\frac{1}{N^{2}}e^{p_{N}}+\varepsilon\biggr)e^{-p_N}\,\drv p_N
\end{equation*}
can be made arbitrarily negative by increasing $b$. Using the upper bound (\ref{equation:exponential-inequality-beta}), it follows that 
\begin{equation*}
\int_{0}^{\infty } \int_{0}^{1} \frac{-xF(\cPriceSDP(x;p_{1},\ldots ,p_{N}))}{1-\beta F(\cPriceSDP(x;p_{1},\ldots ,p_{N}))}\, \drv x \, f(p_{N})\, \drv p_{N}
\end{equation*}
can be made arbitrarily negative by increasing $\beta$ towards $1$, regardless of the values of $p_{1},\ldots ,p_{N-1}$. We therefore conclude that (\ref{equation:exponential-sdp-expected-oos-performance-negative-part}) is unbounded below as $\beta\to 1$.

Recall that $\cPriceSDP(x;p_{1},\ldots ,p_{N})\leq \max\{p_1,\ldots,p_N\}$. The positive term in (\ref{equation:exponential-sdp-expected-oos-performance}) is
\begin{align}\label{equation:exponential-sdp-expected-oos-performance-positive-part}
\Expt_{\Prob^N}\Biggl[\int_{0}^{1}{\frac{\int_{\cPriceSDP(x;\rp_{1},\ldots ,\rp_{N})}^{\infty }pf(p)\, \drv p}{1-\beta F(\cPriceSDP(x;\rp_{1},\ldots ,\rp_{N}))}}\, \drv x\Biggr] &\leq \Expt_{\Prob^N}\Biggl[ \int_{0}^{1}\frac{\int_{\cPriceSDP(x;\rp_{1},\ldots ,\rp_{N})}^{\infty }pf(p)\, \drv p}{1- F(\cPriceSDP(x;\rp_{1},\ldots ,\rp_{N}))}\,  \drv x \Biggr]\\
&= \Expt_{\Prob^N}\Biggl[ \int_{0}^{1} 1+ \cPriceSDP(x;\rp_{1},\ldots ,\rp_{N})  \, \drv x \Biggr]\notag\\
&\leq \Expt_{\Prob^N}\Biggl[ \int_{0}^{1} 1+ \max\{\rp_1,\ldots,\rp_N\} \, \drv x \Biggr]\notag\\
&= \Expt_{\Prob^N}\Bigl[  1+ \max\{\rp_1,\ldots,\rp_N\} \Bigr].\notag
\end{align}
Here $\Expt_{\Prob^N}\bigl[ \max\{\rp_1,\ldots,\rp_N\} \bigr]$ is the expected value of the $N$\textsuperscript{th} order statistic of $N$ random samples from the $\text{Exponential(1)}$ distribution, which is finite; see, e.g., \citep[Section~4.2]{order-statistics:david}. Hence, (\ref{equation:exponential-sdp-expected-oos-performance-positive-part}) is bounded above for all $\beta \in (0,1)$. We have already shown that (\ref{equation:exponential-sdp-expected-oos-performance-negative-part}) is unbounded below as $\beta\to 1$. Thus, the sum (\ref{equation:exponential-sdp-expected-oos-performance}) is unbounded below as $\beta\to 1$, regardless of the choice of $N$.

Consider now $\Expt_{\Prob^N}\bigl[\oExptCtgF_{\MPC}(1; \rp_1,\ldots,\rp_N)\bigr]$. For simplicity, we set $N=2$, although the result can be shown to hold for $N \geq 2$. As for (\ref{equation:exponential-sdp-expected-oos-performance}),
\begin{equation}\label{equation:exponential-mpc-expected-oos-performance}
\Expt_{\Prob^2}\bigl[\oExptCtgF_{\MPC}(1;\rp_1,\rp_2)\bigr] = \Expt_{\Prob^2}\Biggl[\int_{0}^{1}{\frac{\int_{\cPriceMPC(x;\rp_1,\rp_2)}^{\infty }pf(p)\,\drv p-xF(\cPriceMPC(x;\rp_1,\rp_2))}{1-\beta F(\cPriceMPC(x;\rp_1,\rp_2))}}\,\drv x\Biggr],
\end{equation}
which has negative term
\begin{equation}\label{equation:exponential-mpc-expected-oos-performance-negative-part}
\Expt_{\Prob^2}\Biggl[\int_{0}^{1}{\frac{-xF(\cPriceMPC(x;\rp_1,\rp_2))}{1-\beta F(\cPriceMPC(x;\rp_1,\rp_2))}}\,\drv x\Biggr]  = \int_{0}^{\infty } \int_{0}^{\infty }\int_{0}^{1}\frac{-xF(\cPriceMPC(x;p_1,p_2))}{1-\beta F(\cPriceMPC(x;p_1,p_2))}\,\drv x\,f(p_{2})\,\drv p_{2} \,f(p_{1})\,\drv p_{1}.
\end{equation}
The iterated integral (\ref{equation:exponential-mpc-expected-oos-performance-negative-part}) can be divided into ranges based on the value of the sample average $\mu_2 = \frac{1}{2}(p_{1}+p_{2})$. Following (\ref{equation:mpc-policy}), depending on the value of $\mu_2$, either $\cPriceMPC(x;p_1,p_2)=\beta \mu_2-x$ or $\cPriceMPC(x;p_1,p_2)={-x}/{(1-\beta)}$. Since $\mu_2-x \geq \beta \mu_2-x \geq {-x}/{(1-\beta)}$,
\begin{align*}
\int_{0}^{1}\frac{-xF(\cPriceMPC(x;p_1,p_2))}{1-\beta F(\cPriceMPC(x;p_1,p_2))}\,\drv x &\geq \int_{0}^{1}\frac{-xF(\mu_2-x)}{1-F(\mu_2-x)}\,\drv x\\
&= \int_{0}^{\min\{\mu_2,1\}}\frac{-x( 1-e^{-(\mu_2-x)}) }{1-( 1-e^{-(\mu_2-x)}) }\,\drv x \\
&= \frac{1}{2}({\min\{\mu_2,1\}})^2+(1+\min\{\mu_2,1\})e^{\mu_2-{\min\{\mu_2,1\}}}-e^{\mu_2} \\
& \geq -e^{\mu_2}.
\end{align*}
Thus, (\ref{equation:exponential-mpc-expected-oos-performance-negative-part}) is bounded below by
\begin{equation*}
\int_{0}^{\infty } \int_{0}^{\infty }-e^{\frac{1}{2}(p_1+p_2)} e^{-p_{2}}\,\drv p_{2}\,e^{-p_{1}}\,\drv p_{1}=-4
\end{equation*}
for all $\beta\in(0,1)$. Similar reasoning as for (\ref{equation:exponential-sdp-expected-oos-performance-positive-part}) shows that the positive term in (\ref{equation:exponential-mpc-expected-oos-performance}) is bounded above. We therefore conclude that the sum (\ref{equation:exponential-mpc-expected-oos-performance}) is bounded as $\beta\to 1$.
\Halmos
\endproof

\end{APPENDICES}

\ACKNOWLEDGMENT{The authors would like to thank Andrew Mason for helpful discussions. The first and second authors acknowledge support from UOCX2117 MBIE Catalyst Fund New Zealand--German Platform for Green Hydrogen Integration (HINT).}


\bibliographystyle{informs2014} 
\bibliography{bib} 



\end{document}